\documentclass{article}
\usepackage[whole]{bxcjkjatype}

\usepackage{graphicx, amsmath, amsfonts, amssymb, xcolor, amsthm, array, float, subcaption, enumitem}

\usepackage[subrefformat=parens]{subcaption}
\usepackage[margin=2.5cm]{geometry}
\usepackage{hyperref}
\usepackage{mathtools,cleveref}
\usepackage{booktabs,multirow}
\graphicspath{{./fig/}}

\crefname{problem}{Problem}{problem}

\newtheorem{definition}{Definition}[section]

\newtheorem{lemma}[definition]{Lemma}
\newtheorem{theorem}[definition]{Theorem}

\newtheorem{remark}[definition]{Remark}
\newtheorem{example}[definition]{Example}
\newtheorem{conjecture}[definition]{Conjecture}

\newcommand{\R}{\mathbb{R}}

\title{Polyhedra of Constant Gaussian Curvature}
\author{Soto Hisakawa\thanks{Sumitomo Electric, Shimaya, Konohana Ward, Osaka, Japan.
E-mail: \href{mailto:hisakawa.souto@gmail.com}{hisakawa.souto@gmail.com}}
\and
Shizuo Kaji\thanks{Graduate School of Science, Kyoto University. Kitashirakawa Oiwake-cho, Kyoto, Japan.
E-mail: \href{mailto:kaji.shizuo.7r@kyoto-u.ac.jp}{kaji.shizuo.7r@kyoto-u.ac.jp}}
\and
Ryo Kawai\thanks{Fujitsu Ltd.,
Kawasaki, Kanagawa 211-8588, Japan.
E-mail: \href{mailto:kawai.ryo.193@gmail.com}{kawai.ryo.193@gmail.com}}
}
\date{\today}

\begin{document}

\maketitle

\begin{abstract}
Topology and geometry are deeply intertwined in the study of surfaces, though their interaction manifests differently in smooth and discrete settings. In the smooth category, a classical result asserts that any closed smooth surface embedded in $\mathbb{R}^3$ with constant Gaussian curvature must be a sphere, reflecting the strong rigidity of differential geometry.  
In contrast, the discrete setting—where curvature is represented as an angular defect concentrated at vertices—admits far greater flexibility. For instance, a flat torus can be realised as a polyhedral surface in $\mathbb{R}^3$ with zero curvature at every vertex.

We establish a general result: any closed surface, whether orientable or non-orientable and of arbitrary genus, can be realised in $\mathbb{R}^3$ as a (possibly self-intersecting) polyhedral surface in which every vertex has the same angular defect. This highlights a fundamental distinction between discrete and smooth settings, showing that curvature constraints in the discrete realm impose fewer restrictions.
Our proof is constructive and, once recognised, entirely elementary. Yet this fundamental fact appears to have gone unnoticed in the existing literature.
\end{abstract}

\section{Introduction}
It is a classical fact that any $C^2$ isometric embedding of a closed surface into $\mathbb{R}^3$ must have at least one point of positive Gaussian curvature. In particular, although flat metrics exist on the torus (an orientable surface of genus one with $K \equiv 0$) and constant negative curvature metrics exist on closed hyperbolic surfaces (of genus $g \ge 2$), neither can be realised as a $C^2$ surface in $\mathbb{R}^3$. Indeed, Hilbert's theorem~\cite{hilbert1901} shows that no complete surface with constant negative curvature admits an isometric immersion in $\mathbb{R}^3$.

The situation changes when the smoothness requirement is relaxed. The Nash--Kuiper theorem~\cite{kuiper1955c1} asserts that every orientable Riemannian surface admits a $C^1$ isometric embedding into $\mathbb{R}^3$. 
Indeed, $C^1$ flat tori are actually realised in $\mathbb{R}^3$~\cite{flattori}.
Furthermore, 
Burago--Zalgaller's theorem~\cite{burago1996isometric} states that every orientable polyhedral surface that is obtained by gluing a flat polygon boundary has an isometric piecewise linear embedding into $\R^3$.
In particular, explicit
origami realisations of flat tori in $\mathbb{R}^3$ have been constructed~\cite{quintanar2020explicit,tsuboi2024,Lazarus2024}.

A natural way to extend the notion of constant curvature surfaces to the discrete setting is through the Gauss--Bonnet theorem:
\[
\int_S K\, dA = 2\pi \chi(S),
\]
where $dA$ is the area element and $\chi(S)$ denotes the Euler characteristic of the surface $S$. In the polyhedral setting, this relation takes the form of Descartes' theorem~\cite{hilton1982}:
\begin{equation}\label{eq:Descartes}
\sum_{v \in V} \delta_v = 2\pi\chi(P),    
\end{equation}
with the angular defect at a vertex $v$ defined by
\[
\delta_v = 2\pi - \sum \text{(angles incident at } v\text{)}.
\]
This correspondence makes angular defect a natural discrete analogue of Gaussian curvature.
% Some concrete examples of polyhedra with negative angular defects at all vertices are given~\cite{barros2011counterexamples,chern2018shape}.

In this work, we focus on polyhedra with a constant angular defect at every vertex---which we refer to as \emph{constant curvature polyhedra} (CCP). Known constructions include the well-understood cases for orientable surfaces: genus zero (the sphere, as exemplified by the Platonic solids) and genus one (the flat tori discussed above). In addition, there exist several sporadic examples among non-orientable surfaces with genus 1,4,6, and so on. 
Here, we systematically construct both orientable and non-orientable constant curvature polyhedra, thereby realising all topological types for closed surfaces.
Interactive 3D models of constructed polyhedra in GeoGebra and in a 3D printable format are available
\footnote{
\url{https://www.geogebra.org/u/hisakawa}
and
\url{https://github.com/SotoHisakawa/Polyhedra-of-Constant-Gaussian-Curvature/}},
providing visual aids to enhance understanding of the constructions.

%%%%

\section{Polyhedra with self-intersection}

Informally, a polyhedron in $\mathbb{R}^3$ is formed by gluing planar polygons along their edges, with each edge shared by exactly two polygons.
Since we consider both orientable and non-orientable polyhedra of arbitrary genus, we begin by formalising our notion of a polyhedron.

\begin{definition}\label{dfn:polyhedron}
An \emph{abstract polyhedron} $K$ is a finite regular CW complex of dimension~$2$ in which every interior point of a $1$-cell lies in the boundary of exactly two $2$-cells. Regularity means that each attaching map is a homeomorphism from the boundary of the cell to its image in the lower-dimensional skeleton. In particular, the boundary of each $2$-cell is homeomorphic to the circle $S^1$ and attaches along a cyclic sequence of $1$-cells.

We write $K = (V, E, F)$ for the vertex-, edge-, and face-sets (the $0$-, $1$-, and $2$-cells respectively).

A \emph{polyhedron} $P$ in $\mathbb{R}^3$ is the image of an abstract polyhedron $K$ under a continuous map $g: K \to \mathbb{R}^3$ satisfying:
\begin{itemize}
\item For each edge $e \in E$, the image $g(e)$ is a straight line segment joining the images of its endpoints.
\item For each face $f \in F$, the image $g(f)$ is a planar polygon with boundary $g(\partial f)$.
\item The restriction of $g$ to each cell is an embedding.
\item The dihedral angle at each edge is not equal to $\pi$.
\end{itemize}
\end{definition}

The map $g$ is not required to be a global embedding: points belonging to different cells may coincide in $\mathbb{R}^3$, so faces may intersect or pass through one another. A polyhedron is said to be \emph{embedded} if no such global self-intersection occurs.

When the context is clear, we identify each cell of $K$ with its image under $g$.

\begin{remark}
A polyhedron 
may contain singularities such as cross-caps, as in the tetrahemihexahedron (\Cref{fig:n_g1_v6}).
\end{remark}

\medskip

For two polyhedra $P_1 = (V_1, E_1, F_1)$ and $P_2 = (V_2, E_2, F_2)$ with an isometry $\phi$ between a face $f_1 \in F_1$ and a face $f_2 \in F_2$,
we denote their \emph{connected sum along $\phi$} by $P_1 \#_\phi P_2$. 
It satisfies the following:
\begin{enumerate}
    \item $\chi(P_1 \#_\phi P_2) = \chi(P_1) + \chi(P_2) - 2$.
    \item $P_1 \#_\phi P_2$ is orientable if and only if both $P_1$ and $P_2$ are orientable.
\end{enumerate}

We recall the classification of connected closed surfaces.

\begin{lemma}\label{lem:classification}
Every connected closed surface $S$ is homeomorphic to exactly one of:
\begin{enumerate}
\item the sphere $S^2$;
\item a connected sum of $g$ tori $T^2$, where $g \ge 1$;
\item a connected sum of $g$ projective planes $\mathbb{RP}^2$, where $g \ge 1$.
\end{enumerate}
In both the orientable and non-orientable cases, $g$ is called the \emph{genus} of~$S$.

Moreover, the Euler characteristic $\chi$ and orientability determine the topological type:
\begin{itemize}
\item Sphere: $\chi = 2$;
\item Orientable genus $g$: $\chi = 2 - 2g$;
\item Non-orientable genus $g$: $\chi = 2 - g$.
\end{itemize}
\end{lemma}

\medskip

We now introduce an operation on a polyhedron that realises the topological attachment of a handle.

\begin{definition}[Drilling]
Let $P$ be a polyhedron containing two parallel faces $f_1$ and $f_2$. Let $\ell$ be a line segment orthogonal to both $f_1$ and $f_2$ and intersecting $P$:
\begin{itemize}
\item only in the interiors of $f_1$ and $f_2$;
\item exactly once in each, at points $p_1 \in f_1$ and $p_2 \in f_2$.
\end{itemize}

For an integer $n \ge 3$ and sufficiently small $\varepsilon > 0$, the \emph{drilling operation} along $\ell$ consists of:
\begin{enumerate}
\item Constructing a regular $n$-gonal prism $\Pi(\ell, n)$ of radius $\varepsilon$ whose axis is $\ell$, with bases centred at $p_1$ and $p_2$ lying in $f_1$ and $f_2$ respectively.
\item Removing from each of $f_1$ and $f_2$ the interior of the corresponding $n$-gon.
\item Attaching the side faces of $\Pi(\ell, n)$ to $P$, thereby forming a tunnel between $f_1$ and $f_2$.
\item Subdividing the resulting $n$-gons with $n$-gonal holes into polygonal 2-cells using only existing vertices. The subdivision is not canonical, and any choice suffices.
\end{enumerate}

The resulting polyhedron, denoted $P \# \Pi(\ell, n)$, has Euler characteristic $\chi(P) - 2$ and $2n$ additional vertices, each with angular defect $\delta' = -\tfrac{2\pi}{n}$. The final subdivision and the choice of $\varepsilon$ are immaterial to the construction.
\end{definition}

We now describe the main construction for constant-curvature polyhedra (CCPs). Let $K = (V, E, F)$ be an abstract polyhedron. By Descartes' theorem~\eqref{eq:Descartes}, a polyhedron with $V$ vertices and Euler characteristic $\chi$ must have constant defect
\[
\delta = \frac{2\pi\,\chi}{V} = \frac{2\pi\,(V - E + F)}{V}.
\]
For small genus $g$, we construct explicit realisations by assigning coordinates to vertices in $\mathbb{R}^3$ so that each vertex achieves defect $\delta$, ensuring that the resulting CCP is convenient for subsequent modifications.

Suppose $P \subset \mathbb{R}^3$ is a genus-$g$ CCP with $V$ vertices and that it contains a pair of parallel faces allowing drilling along some segment~$\ell$. For integers $n \ge 3$ and $k \ge 1$, where $k$ is the number of drills, set
\[
P' = P \# \Pi(\ell_1, n) \# \cdots \# \Pi(\ell_k, n) \quad (\text{$k$ times}),
\]
where the $\ell_i$ are parallel copies of~$\ell$ contained in a small neighbourhood of it. Then,
 $V' = V + 2nk$ and
 $\chi(P') = \chi(P) - 2k$.

Consequently:
\begin{itemize}
\item If $P$ is orientable of genus $g$, then $P'$ is orientable of genus $g + k$.
\item If $P$ is non-orientable of genus $g$, then $P'$ is non-orientable of genus $g + 2k$.
\end{itemize}

If the integer
\begin{equation}\label{eq:choice-of-n}
n = -\frac{V}{\chi(P)}
\end{equation}
exists, then $P'$ again has constant defect $\delta$ and is therefore a CCP of the same defect and higher genus.

\medskip

The paper is organised into two parts. First, we establish the existence of CCPs for all closed surfaces (\Cref{thm:main}). Second, we investigate CCPs with minimal vertex counts and propose conjecturally optimal constructions (\Cref{thm:minimal}).

%%%%%%%%%%%%%%%%%%
\section{Construction of Polyhedra}

Our first primary result is the construction of a CCP
for every topological type.
\begin{theorem}\label{thm:main}
Constant curvature polyhedra exist for every genus, regardless of orientability.

Moreover, the orientable ones are realised without self-intersection.
\end{theorem}

Based on the classification of surfaces (\Cref{lem:classification}), we give explicit constructions for every topological type.

\subsection{Orientable case}

For genus $g=0$ (the sphere), the regular tetrahedron is a CCP.

\begin{figure}[htbp]
    \centering
    \includegraphics[width=0.2\linewidth]{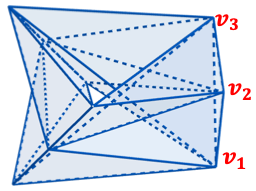} \hfill
    \includegraphics[width=0.2\linewidth]{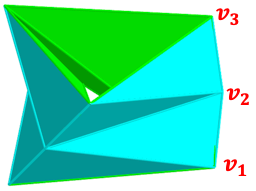} \hfill
    \includegraphics[width=0.2\linewidth]{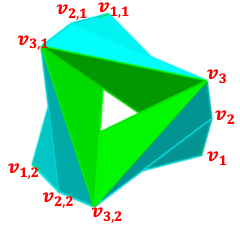}\\
    \includegraphics[width=0.5\linewidth]{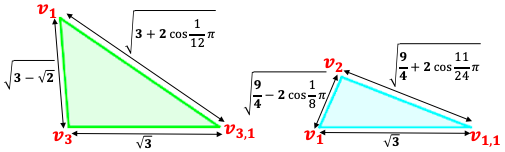}\hfill
    \includegraphics[width=0.4\linewidth]{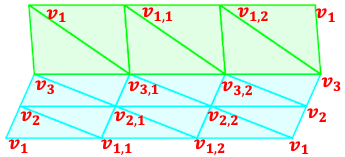}
    \caption{
        Top: A flat torus with nine vertices.
        Bottom: A net of the flat torus.
        At each vertex, six triangles meet, forming a tiling of the plane with zero angular defect.  
        Specifically, there are two types of triangles, with three of each type meeting at every vertex so that the angles sum to $2\pi$.
    }
    \label{fig:o_g1}
\end{figure}

For genus $g=1$ (the torus), the origami tori constructed in \cite{tsuboi2024} are CCP.
For the sake of completeness, we detail a construction of a flat torus with 9 vertices
(\Cref{fig:o_g1}).
It has three-fold rotational symmetry around the $z$-axis.
The coordinates of the vertices are given as
\[
\begin{aligned}
v_1 &= (1, 0, 0), \\
v_2 &= \left(\cos{\frac{1}{8}\pi},\ \sin{\frac{1}{8}\pi},\ \frac{1}{2}\right), \\
v_3 &= \left(\cos{\frac{1}{4}\pi},\ \sin{\frac{1}{4}\pi},\ 1\right),
\end{aligned}
\]
and their images under rotations by $2\pi/3$ and
      $4\pi/3$ about the z-axis, denoted $v_{1,1}, v_{2,1}, v_{3,1}$ and $v_{1,2}, v_{2,2}, v_{3,2}$ respectively.
      
The lengths of the edges are computed as follows:
\begin{align*}
| v_1 v_3 | &= \left| (1,0,0),\ \left(\cos{\frac{1}{4}\pi},\ \sin{\frac{1}{4}\pi},\ 1\right) \right|
= \sqrt{3 - \sqrt{2}} \\[1em]
| v_1 v_{3,1} | &= \left| (1,0,0),\ \left(\cos{\frac{11}{12}\pi},\ \sin{\frac{11}{12}\pi},\ 1\right) \right|
= \sqrt{3 + 2\cos{\frac{1}{12}\pi}} \\[1em]
| v_3 v_{3,1} | &= \left| \left(\cos{\frac{1}{4}\pi},\ \sin{\frac{1}{4}\pi},\ 1\right),\ \left(\cos{\frac{11}{12}\pi},\ \sin{\frac{11}{12}\pi},\ 1\right) \right|
=\sqrt{3} \\[1em]
| v_2 v_1 | &= \left| \left(\cos{\frac{1}{8}\pi},\ \sin{\frac{1}{8}\pi},\ \frac{1}{2}\right),\ (1,0,0) \right|
= \sqrt{\frac{9}{4} - 2\cos{\frac{1}{8}\pi}} \\[1em]
| v_2 v_{1,1} | &= \left| \left(\cos{\frac{1}{8}\pi},\ \sin{\frac{1}{8}\pi},\ \frac{1}{2}\right),\ \left(-\frac{1}{2},\ \frac{\sqrt{3}}{2},\ 0\right) \right|
= \sqrt{\frac{9}{4} + 2\cos{\frac{11}{24}\pi}} \\[1em]
| v_1v_{1,1} | &= \left| (1,0,0),\ \left(-\frac{1}{2},\ \frac{\sqrt{3}}{2},\ 0\right) \right|
= \sqrt{3}.
\end{align*}

At each vertex, three copies of each of the two triangle types or six copies of a triangle meet, and the sum of their angles is $2\pi$,
thereby guaranteeing zero angular defect (flatness) at every vertex.

\begin{figure}[htbp]
    \centering
    \includegraphics[width=0.2\linewidth]{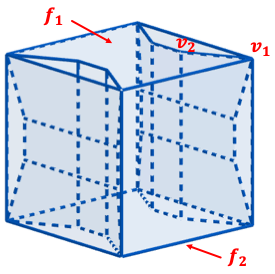}\hfill
    \includegraphics[width=0.2\linewidth]{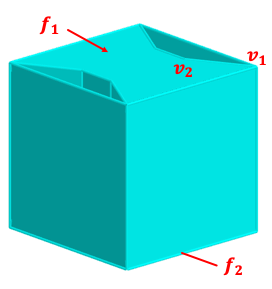}\hfill
    \includegraphics[width=0.2\linewidth]{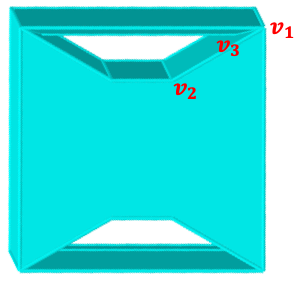}\hfill
    \includegraphics[width=0.2\linewidth]{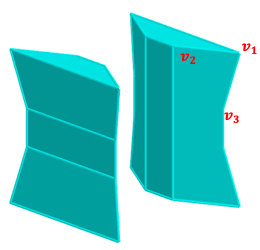} \\
    \includegraphics[width=0.8\linewidth]{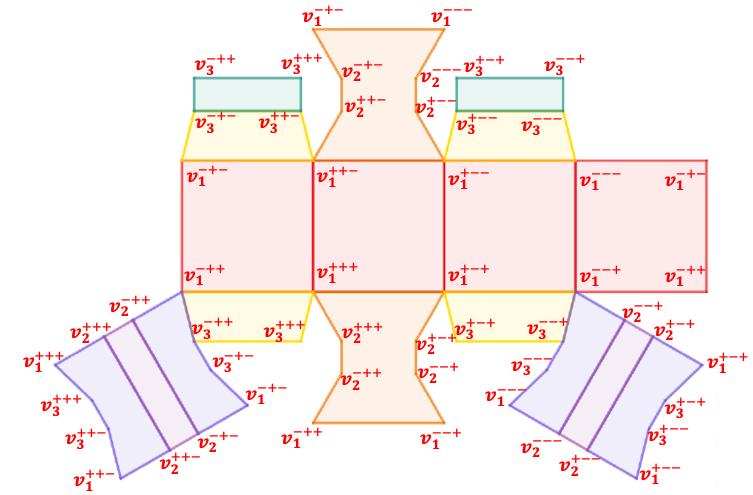}\hfill
    \includegraphics[width=0.9\linewidth]{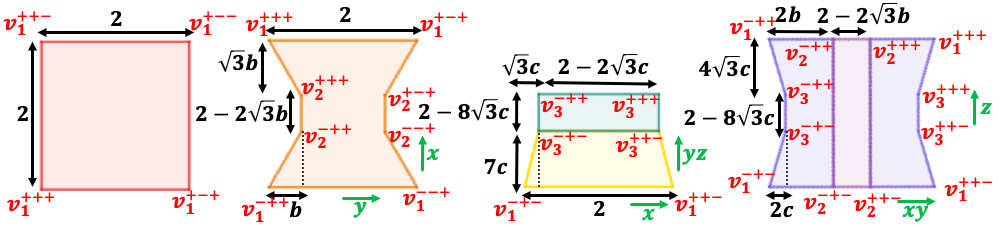}
    \caption{Top: An orientable genus-$2$ CCP $P^2_{24}$ with 24 vertices and its two holes.
    Bottom: A net of $P^2_{24}$.}
    \label{fig:o_g2_v24_k}
\end{figure}

\medskip
Next, we construct an orientable CCP of genus \(g = 2\) with 24 vertices, denoted by \(P^2_{24}\), shown in \Cref{fig:o_g2_v24_k} (Top).  
We begin by outlining the idea.

To simplify the computation of the angular defect \(\delta_v\) at each vertex, we consider a polyhedron with high symmetry, namely one symmetric with respect to the coordinate planes \(xy\), \(yz\), and \(zx\).  
As our base polyhedron we take the cube, and then introduce two holes so that the resulting genus is \(g=2\).  

Let \(v_1\) be a cube vertex in the first octant, \(v_2\) a point on a cube face, and \(v_3\) a point inside the cube, with the combinatorial structure indicated in \Cref{fig:o_g2_v24_k} (Top left).  
By \eqref{eq:Descartes}, the angular defect at each vertex must be
\[
\delta_v = \frac{2\pi(2-2g)}{|V|} = -\frac{\pi}{6}.
\]

We first place \(v_2\) so that \(\delta_{v_2} = -\frac{\pi}{6}\), and then adjust the location of \(v_3\) accordingly.  
More precisely, let \(b, c \in \mathbb{R}\) satisfy
\[
b > c > 0, 
\qquad 1 > \sqrt{3}b,
\qquad 1 > 4\sqrt{3}c.
\]
We then define the three vertices in the first octant as
\[
v_1 = (1,\, 1,\, 1), \quad
v_2 = (1 - \sqrt{3}b,\, 1 - b,\, 1), \quad
v_3 = (1 - \sqrt{3}c,\, 1 - c,\, 1 - 4\sqrt{3}c).
\]
We denote by $v^{\sigma_1\sigma_2\sigma_3}$ the image of a vertex $v$ under reflection with respect to the coordinate planes, where $\sigma_i\in\{+,-\}$ indicates the sign of the $i$-th coordinate. For example, $v_1^{++-}=(1,1,-1)$. Accordingly, we write the symmetric copies of $v_1,v_2,v_3$ as $v_j^{\sigma_1\sigma_2\sigma_3}$.
%Their symmetric copies are denoted, for instance, by \(v_1^{++-}\), where the signs (\(\mathrm{++-}\)) indicate the octant containing the vertex.
%Accordingly, we write the symmetric copies of $v_1,v_2,v_3$ as $v_j^{\sigma_1\sigma_2\sigma_3}$.

% Since \(b > c\), the edges \(v_2^{+++}v_2^{++-}\) and \(v_3^{+++}v_3^{-++}\) do not intersect, producing a hole as shown in \Cref{fig:o_g2_v24_k} (Right). 
Consequently the resulting polyhedron has two holes, as shown in \Cref{fig:o_g2_v24_k} (Right).

Counting vertices, edges, and faces gives
\[
|V| = 24, \quad |E| = 44, \quad |F| = 18,
\]
and hence
\[
\chi(P^2_{24}) = |V| - |E| + |F| = -2 = 2 - 2g,
\]
confirming that \(g=2\).

Finally, we verify \(\delta_v = -\frac{\pi}{6}\) by explicitly computing edge lengths and applying the law of cosines.  
For clarity, we use the notation \(|uv|_{yz}\) to denote the distance between vertices \(u\) and \(v\) after projection onto the \(yz\)-plane.  

\begingroup\setlength{\jot}{.25ex}
\[
\begin{alignedat}{2}
|v_1^{+++}v_2^{+++}|_{y} &= \bigl| 1,\ 1-b \bigr| = b, &\quad
|v_1^{+++}v_2^{+++}|_{x} &= \bigl| 1,\ 1-\sqrt{3}b \bigr| = \sqrt{3}b, \\
|v_1^{+++}v_3^{+++}|_{z} &= \bigl| 1,\ 1-4\sqrt{3}c \bigr| = 4\sqrt{3}c, &\quad
|v_1^{+++}v_3^{+++}|_{xy} &= \bigl| (1,1),\ (1-\sqrt{3}c,1-c) \bigr| = 2c, \\
|v_1^{+++}v_3^{+++}|_{yz} &= \bigl| (1,1),\ (1-c,1-4\sqrt{3}c) \bigr| = 7c, &\quad
|v_1^{+++}v_3^{+++}|_{x} &= \bigl| 1,\ 1-\sqrt{3}c \bigr| = \sqrt{3}c.
\end{alignedat}
\]
\endgroup

% \[
% \angle v_1^{+-+}v_1^{+++}v_2^{+++} = \arctan{\frac{\|v_1^{+++}v_2^{+++}\|_x}{\|v_1^{+++}v_2^{+++}\|_y} = \arctan{\frac{\sqrt{3}b}{b}}} = \frac{\pi}{3}
% \]
% \[
% \angle v_2^{+++}v_1^{+++}v_3^{+++} = \frac{\|v_1^{+++}v_3^{+++}\|_{z}}{\|v_1^{+++}v_3^{+++}\|_{xy}} = \arctan{\frac{4\sqrt{3}c}{2c}} = \arctan{\left(2\sqrt{3}\right)}
% \]
% \[
% \angle v_3^{+++}v_1^{+++}v_1^{-++} = \arctan{\frac{\|v_1^{+++}v_3^{+++}\|_{yz}}{\|v_1^{+++}v_3^{+++}\|_{x}}} = \arctan{\frac{7c}{\sqrt{3}c}} = \arctan{\frac{7}{\sqrt{3}}}
% \]
\begingroup\setlength{\jot}{.25ex}
\[
\begin{aligned}
\angle v_2^{-++}v_2^{+++}v_1^{+++}
  &= \pi + \arctan\frac{|v_1^{+++}v_2^{+++}|_{y}}{|v_1^{+++}v_2^{+++}|_{x}}
   = \pi + \arctan\frac{b}{\sqrt{3}b}
   = \frac{7\pi}{6},\\
\angle v_1^{+++}v_3^{+++}v_3^{++-}
  &= \pi + \arctan\frac{|v_1^{+++}v_3^{+++}|_{xy}}{|v_1^{+++}v_3^{+++}|_{z}}
   = \pi + \arctan\frac{2c}{4\sqrt{3}c}
   = \pi + \arctan\frac{1}{2\sqrt{3}},\\
\angle v_3^{-++}v_3^{+++}v_1^{+++}
  &= \frac{\pi}{2} + \arctan\frac{|v_1^{+++}v_3^{+++}|_{x}}{|v_1^{+++}v_3^{+++}|_{yz}}
   = \frac{\pi}{2} + \arctan\frac{\sqrt{3}c}{7c}
   = \frac{\pi}{2} + \arctan\frac{\sqrt{3}}{7} .
\end{aligned}
\]
\endgroup

Based on above, the angular defect at the vertices can be calculated as follows:

% \begin{align*}
% \delta_{v_1} &= 2\pi - \bigg(
% \angle v_1^{-++}v_1^{+++}v_1^{++-}
% + \angle v_1^{++-}v_1^{+++}v_1^{+-+}
% + \angle v_1^{+-+}v_1^{+++}v_2^{+++} \\
% &\hspace{4em} + \angle v_2^{+++}v_1^{+++}v_3^{+++}
% + \angle v_3^{+++}v_1^{+++}v_1^{-++}
% \bigg) \\
% &= 2\pi - \left(
% \frac{\pi}{2} + \frac{\pi}{2}
% + \arctan\frac{\sqrt{3}b}{b}
% + \arctan\frac{7c}{\sqrt{3}c}
% + \arctan\frac{4\sqrt{3}c}{2c}
% \right)
% \end{align*}

\begin{align*}
\delta_{v_2} &= 2\pi - \left(
\angle v_1^{+++}v_2^{+++}v_2^{++-}
+ \angle v_2^{++-}v_2^{+++}v_2^{-++}
+ \angle v_2^{-++}v_2^{+++}v_1^{+++}
\right) \\
&= 2\pi - \left(
\frac{\pi}{2} + \frac{\pi}{2} + \frac{7\pi}{6}
\right)\\
&= -\frac{\pi}{6}, \\[1em]
\delta_{v_3} &= 2\pi - \left(
\angle v_1^{+++}v_3^{+++}v_3^{++-}
+ \angle v_3^{++-}v_3^{+++}v_3^{-++}
+ \angle v_3^{-++}v_3^{+++}v_1^{+++}
\right) \\
&= 2\pi - \left(
\left(\pi + \arctan\frac{1}{2\sqrt{3}}\right)
+ \frac{\pi}{2}
+ \left(\frac{\pi}{2} + \arctan{\frac{\sqrt{3}}{7}}\right)
\right)\\
&= -\frac{\pi}{6},
\end{align*}
% Using $\arctan\sqrt{3} = \frac{\pi}{3}$ and $\arctan\frac{1}{\sqrt{3}} = \frac{\pi}{6}$, and 
where we used
\[
\arctan\left(\frac{1}{2\sqrt{3}}\right) + \arctan\left(\frac{\sqrt{3}}{7}\right)
= \arctan\left(\frac{1}{\sqrt{3}}\right)
= \frac{\pi}{6}
\]
by the addition theorem for arctangent:
\[
\arctan\alpha + \arctan\beta = \arctan\left( \frac{\alpha + \beta}{1 - \alpha\beta} \right).
\]
Since we verified \(\delta_{v_2} = \delta_{v_3} = -\pi/6\), and by Descartes' theorem \eqref{eq:Descartes} the sum of all defects equals \(2\pi\chi(P)\), we conclude \(\delta_{v_1} = -\pi/6\) as well.
Therefore, all vertices have the same angular defect of
\(\delta_{v}=-\frac{\pi}{6}\).

For orientable CCPs of higher genus, the genus can be increased by performing the drilling operation. Specifically, by drilling through a line segment~$\ell$ connecting the upper and lower faces $f_1$ and $f_2$, and inserting a prism $\Pi(\ell, 12)$ (since $12=-V/\chi$ by \eqref{eq:choice-of-n}), we increase the genus by one without self-intersection.  
By repeating this drilling operation $g-2$ times, we obtain an orientable genus-$g$ CCP, denoted by $P^g_{24(g-1)}$.

%%%%%%%%%%%%%%%%%%%%%%%%%%
\subsection{Non-orientable case}

\begin{figure}[htbp]
    \centering
    \includegraphics[width=0.2\linewidth]{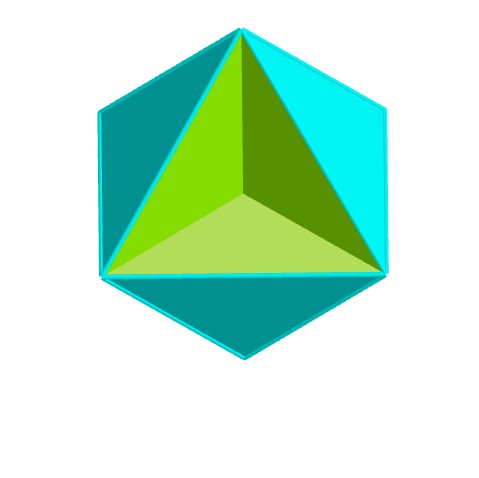}\hfill
    \includegraphics[width=0.2\linewidth]{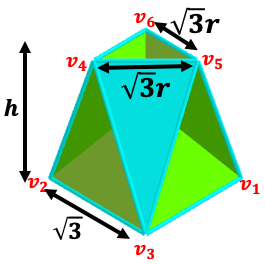} \hfill
    \includegraphics[width=0.3\linewidth]{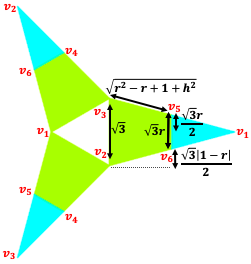}\hfill
    \includegraphics[width=0.3\linewidth]{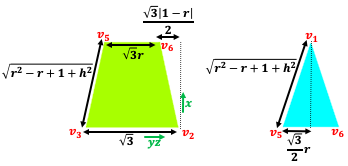}
    \caption{Left: The tetrahemihexahedron $Q^1_6$ with the topology of the real projective plane. Centre and Right: A modified tetrahemihexahedron $R(r,h)$.}
    \label{fig:n_g1_v6}
\end{figure}

For genus $g=1$ (the real projective plane), the tetrahemihexahedron $Q^1_6$ with 
$|V| = 6$, $|E| = 12$, and $|F| = 7$, as shown in \Cref{fig:n_g1_v6} (Left), is a CCP.
In the tetrahemihexahedron, two equilateral triangles and two squares meet at each vertex, resulting in an angular defect of $\frac{\pi}{3}$ at each vertex.

For genus $g=2$ (the Klein bottle), a CCP $Q^2_9$ can be constructed as the connected sum of two projective planes, with the parameters chosen so that the angular defect at each vertex is zero.
To achieve this, we introduce a polyhedron, denoted $R(r,h)$, obtained by modifying the edge lengths of the tetrahemihexahedron as shown in \Cref{fig:n_g1_v6} (Centre and Right).
We set $v_1 = (1, 0, 0)$ and $v_6 = \left(\frac{1}{2} r,\, \frac{\sqrt{3}}{2} r,\, h\right)$, and obtain the other vertices by rotating these points by angles $\frac{\pi}{3}$ and $\frac{2\pi}{3}$ about the $z$-axis.
Each quadrilateral face is an isosceles trapezoid with a top base of $\sqrt{3} r$, a bottom base of $\sqrt{3}$, and a height of $\sqrt{\frac{1}{4}(r+1)^2 + h^2}$.
Each triangular face is an isosceles triangle with a base of $\sqrt{3} r$ and a height of $\sqrt{\frac{1}{4} r^2 - r + 1 + h^2}$.
We have the following edge lengths:
\begin{align*}
|v_1 v_5| &= \left| (1, 0, 0),\, \left( \frac{1}{2} r,\, -\frac{\sqrt{3}}{2} r,\, h \right) \right| = \sqrt{r^2 - r + 1 + h^2}, \\
|v_5 v_6| &= \left| \left( \frac{1}{2} r,\, -\frac{\sqrt{3}}{2} r,\, h \right),\, \left( \frac{1}{2} r,\, \frac{\sqrt{3}}{2} r,\, h \right) \right| = \sqrt{3} r, \\
|v_5 v_3|_{yz} &= \left| \left( \frac{1}{2} r,\, -\frac{\sqrt{3}}{2} r,\, h \right),\, \left( -\frac{1}{2},\, -\frac{\sqrt{3}}{2},\, 0 \right) \right| = \frac{\sqrt{3} |1 - r|}{2}.
\end{align*}
% \begin{align*}
% \|v_3 v_2\| &= \left\| \left( -\frac{1}{2},\, -\frac{\sqrt{3}}{2},\, 0 \right),\, \left( -\frac{1}{2},\, \frac{\sqrt{3}}{2},\, 0 \right) \right) = \sqrt{3}, \\
% \|v_5 v_3\| &= \left\| \left( \frac{1}{2} r,\, -\frac{\sqrt{3}}{2} r,\, h \right),\, \left( -\frac{1}{2},\, -\frac{\sqrt{3}}{2},\, 0 \right) \right) = \sqrt{r^2 - r + 1 + h^2}\\
% \angle v_5 v_1 v_6 &= \arccos{\frac{\|v_1v_5\|^2+\|v_1v_6\|^2-\|v_5v_6\|^2}{2\|v_1v_5\|\|v_1v_6\|}}\\
%     &= \arccos{\frac{\|v_1v_5\|^2+\|v_1v_5\|^2-\|v_5v_6\|^2}{2\|v_1v_5\|\|v_1v_5\|}}\\
%     &= \arccos{\frac{2\|v_1v_5\|^2-\|v_5v_6\|^2}{2\|v_1v_5\|^2}}\\
%     &= \arccos{\frac{2(r^2 - r + 1 + h^2)-(\sqrt{3}r)^2}{2(r^2 - r + 1 + h^2)}}\\
%     &= \arccos{\frac{-r^2 - 2r + 2 + 2h^2}{2(r^2 - r + 1 + h^2)}},\\[1em]
% \angle v_6 v_5 v_3 &= \pi - \arccos{\frac{\|v_6 v_3\|_{yz}}{\|v_6 v_3\|}}\\
%     &= \pi - \arccos{\frac{\|v_6 v_3\|_{yz}}{\|v_1 v_5\|}}\\
%     &= \pi - \arccos{\frac{\frac{\sqrt{3} (1 - r)}{2}}{\sqrt{r^2 - r + 1 + h^2}}}\\
%     &= \pi - \arccos{\frac{\sqrt{3} (1 - r)}{2\sqrt{r^2 - r + 1 + h^2}}},\\[1em]
% \end{align*}
The relevant angles are computed as:
\begin{align*}
\angle v_1 v_5 v_6 &= \arccos{\frac{\frac{|v_5 v_6|}{2}}{|v_1 v_5|}}
%    &= \arccos{\frac{\frac{\sqrt{3}r}{2}}{\sqrt{r^2 - r + 1 + h^2}}}\\
    = \arccos{\frac{\sqrt{3}r}{2\sqrt{r^2 - r + 1 + h^2}}},\\[1em]
\angle v_5 v_3 v_2 &= \arccos{\frac{|v_5 v_3|_{yz}}{|v_5 v_3|}}
    = \arccos{\frac{|v_5 v_3|_{yz}}{|v_1 v_5|}}\\
    %&= \arccos{\frac{\frac{\sqrt{3} (1 - r)}{2}}{\sqrt{r^2 - r + 1 + h^2}}}\\
    &= \arccos{\frac{\sqrt{3} (1 - r)}{2\sqrt{r^2 - r + 1 + h^2}}}
\end{align*}

% The sum of the angles at $v_1$ is given by
% \begin{align*}
% \sum \text{(angles at $v_1$)} &= \angle v_3 v_1 v_6  + \angle v_6 v_1 v_5 + \angle v_5 v_1 v_2 \\
% &= \arccos{\frac{ \frac{\sqrt{3} - \sqrt{3} r }{2} }{ \sqrt{ r^2 - r + 1 + h^2 } } }
%   + \arccos{\frac{ \frac{\sqrt{3} - \sqrt{3} r }{2} }{ \sqrt{ r^2 - r + 1 + h^2 } } } \\
% &\quad + \arccos{\frac{ 2(r^2 - r + 1 + h^2) - (\sqrt{3} r)^2 }{ 2 ( r^2 - r + 1 + h^2 ) } } \\
% &= 2\,\arccos{\frac{ \frac{\sqrt{3} - \sqrt{3} r }{2} }{ \sqrt{ r^2 - r + 1 + h^2 } } } + \arccos{\frac{ 2(r^2 - r + 1 + h^2) - (\sqrt{3} r)^2 }{ 2 ( r^2 - r + 1 + h^2 ) } }.
% \end{align*}
% Similarly, the sum of the angles at $v_6$ is
% \begin{align*}
% \sum \text{(angles at $v_6$)} &= \angle v_1 v_6 v_4 + \angle v_4 v_6 v_2 + \angle v_2 v_6 v_5 + \angle v_5 v_6 v_1 \\
% &= \arccos{ \frac{ \frac{\sqrt{3} r - \sqrt{3} }{2} }{ \sqrt{ r^2 - r + 1 + h^2 } } }
%   + \arccos{ \frac{ \frac{\sqrt{3} r }{2} }{ \sqrt{ r^2 - r + 1 + h^2 } } } \\
% &\quad + \arccos{ \frac{ \frac{\sqrt{3} r - \sqrt{3} }{2} }{ \sqrt{ r^2 - r + 1 + h^2 } } }
%   + \arccos{ \frac{ \frac{\sqrt{3} r }{2} }{ \sqrt{ r^2 - r + 1 + h^2 } } } \\
% &= 2\,\arccos{ \frac{ \frac{\sqrt{3} r - \sqrt{3} }{2} }{ \sqrt{ r^2 - r + 1 + h^2 } } }
%   + 2\,\arccos{ \frac{ \frac{\sqrt{3} r }{2} }{ \sqrt{ r^2 - r + 1 + h^2 } } }.
% \end{align*}

\begin{figure}[htbp]
    \centering
    \includegraphics[width=0.2\linewidth]{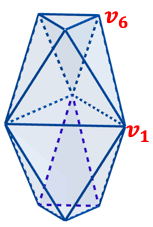}\hfill
    \includegraphics[width=0.2\linewidth]{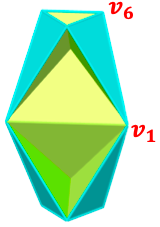}\hfill
    \includegraphics[width=0.2\linewidth]{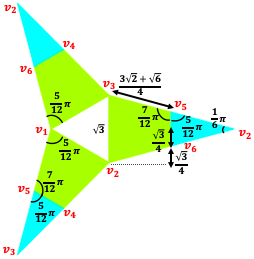}
    \caption{A flat Klein bottle $Q^2_9$.}
    \label{fig:n_g2_v9}
\end{figure}

Taking the connected sum of two copies of 
$R\left(r, h\right)$ with $r = \frac{1}{2}$ and $h = \frac{1}{2}\sqrt{3(1 + \sqrt{3})}$ along the face $v_1 v_2 v_3$ yields the polyhedron 
\[
Q^2_9 = R(r, h) \# R(r, h)
\]
depicted in \Cref{fig:n_g2_v9}, which is a non-orientable genus-$2$ CCP.
It is easy to check that $|V| = 9$, $|E| = 21$, and $|F| = 12$, and thus,
\[
2 - g = \chi(Q^2_9) = |V| - |E| + |F| = 0.
\]
Therefore, the angular defect at each vertex is
\[
\delta_v = \frac{2\pi (2 - g)}{9} = 0,
\]
in accordance with \eqref{eq:Descartes}.
Indeed, we have
\[
\sqrt{r^2 - r + 1 + h^2} = \sqrt{\left(\frac{1}{2}\right)^2 - \left(\frac{1}{2}\right) + 1 + \frac{3(1 + \sqrt{3})}{4}} = \sqrt{\left(\frac{\sqrt{3}(\sqrt{3}+1)}{2}\right)^2} = \frac{\sqrt{3}(\sqrt{3}+1)}{2}
\]
and
\[
\angle v_1 v_5 v_6 = \arccos{\frac{\sqrt{3}r}{2\sqrt{r^2 - r + 1 + h^2}}} = \arccos{\frac{\sqrt{3}\left(\frac{1}{2}\right)}{2\frac{\sqrt{3}(\sqrt{3}+1)}{2}}}=\arccos{\frac{\sqrt{6}-\sqrt{2}}{4}}= \frac{5\pi}{12}
\]
\[
\angle v_5v_1v_6 = \pi - (\angle v_1 v_5 v_6 + \angle v_1 v_6 v_5) = \pi - (\angle v_1 v_5 v_6 + \angle v_1 v_5 v_6) = \pi - 2\,\angle v_1 v_5 v_6 = \pi - 2\,\frac{5\pi}{12} = \frac{\pi}{6}
\]
\[
\angle v_5 v_3 v_2 = \arccos{\frac{\sqrt{3} \left(1 - r\right)}{2\sqrt{r^2 - r + 1 + h^2}}} = \arccos{\frac{\sqrt{3} \left(1 - \frac{1}{2}\right)}{2\frac{\sqrt{3}(\sqrt{3}+1)}{2}}}=\arccos{\frac{\sqrt{6}-\sqrt{2}}{4}}= \frac{5\pi}{12}
\]
\[
\angle v_3 v_5 v_6 = \pi - \angle v_5 v_3 v_2 = \pi - \frac{5\pi}{12}= \frac{7\pi}{12}
\]

% \begin{align*}
% \angle v_1 v_5 v_6 &= \arccos{\frac{\sqrt{3}r}{2\sqrt{r^2 - r + 1 + h^2}}}\\
%     &= \arccos{\frac{\sqrt{3}\left(\frac{1}{2}\right)}{2\sqrt{\left(\frac{1}{2}\right)^2 - \left(\frac{1}{2}\right) + 1 + \frac{3(1 + \sqrt{3})}{4}}}}\\
%     &= \arccos{\frac{\sqrt{3}\left(\frac{1}{2}\right)}{2\sqrt{\left(\frac{\sqrt{3}(\sqrt{3}+1)}{2}\right)^2}}}\\
%     &= \arccos{\frac{\sqrt{3}\left(\frac{1}{2}\right)}{\sqrt{3}(\sqrt{3}+1)}}\\
%     &= \arccos{\frac{\sqrt{6}-\sqrt{2}}{4}}\\
%     &= \frac{5\pi}{12},\\[1em]
% \angle v_5v_1v_6 &= \pi - (\angle v_1 v_5 v_6 + \angle v_1 v_6 v_5)\\
%     &= \pi - (\angle v_1 v_5 v_6 + \angle v_1 v_5 v_6)\\
%     &= \pi - 2\,\angle v_1 v_5 v_6\\
%     &= \pi - 2\,\frac{5\pi}{12}\\
%     &= \frac{\pi}{6}\\
% \angle v_5 v_3 v_2 &= \arccos{\frac{\sqrt{3} \left(1 - r\right)}{2\sqrt{r^2 - r + 1 + h^2}}}\\
%     &= \arccos{\frac{\sqrt{3} \left(1 - \left(\frac{1}{2}\right)\right)}{2\sqrt{\left(\frac{1}{2}\right)^2 - \left(\frac{1}{2}\right) + 1 + \frac{3(1 + \sqrt{3})}{4}}}}\\
%     &= \arccos{\frac{\sqrt{3} \left(\frac{1}{2}\right)}{2\sqrt{\left(\frac{\sqrt{3}(\sqrt{3}+1)}{2}\right)^2}}}\\
%     &= \arccos{\frac{\sqrt{3}\left(\frac{1}{2}\right)}{\sqrt{3}(\sqrt{3}+1)}}\\
%     &= \arccos{\frac{\sqrt{6}-\sqrt{2}}{4}}\\
%     &= \frac{5\pi}{12},\\[1em]
% \angle v_3 v_5 v_6 &= \pi - \angle v_5 v_3 v_2\\
%     &= \pi - \frac{5\pi}{12}\\
%     &= \frac{7\pi}{12}.
% \end{align*}
Therefore, 
\begin{align*}
\delta_{v_1} &= 2\pi - 2\left( \angle v_3 v_1 v_6 + \angle v_6 v_1 v_5 + \angle v_5 v_1 v_2 \right) \\
             &= 2\pi - 2\left( \angle v_5 v_3 v_2 + \angle v_5 v_1 v_6 + \angle v_5 v_3 v_2 \right) \\
             &= 2\pi - 2\left( \frac{5\pi}{12} + \frac{\pi}{6} + \frac{5\pi}{12} \right) = 0, \\
\delta_{v_6} &= 2\pi - \left( \angle v_1 v_6 v_4 + \angle v_4 v_6 v_2 + \angle v_2 v_6 v_5 + \angle v_5 v_6 v_1 \right) \\
             &= 2\pi - \left( \angle v_3 v_5 v_6 + \angle v_1 v_5 v_6 + \angle v_3 v_5 v_6 + \angle v_1 v_5 v_6 \right) \\
             &= 2\pi - \left( \frac{5\pi}{12} + \frac{7\pi}{12} + \frac{5\pi}{12} + \frac{7\pi}{12} \right) = 0.
\end{align*}

\begin{figure}[htbp]
    \centering
    \includegraphics[width=0.2\linewidth]{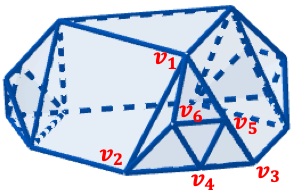} \hfill
    \includegraphics[width=0.2\linewidth]{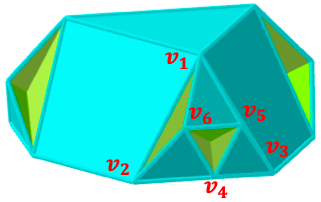}\hfill
    \includegraphics[width=0.2\linewidth]{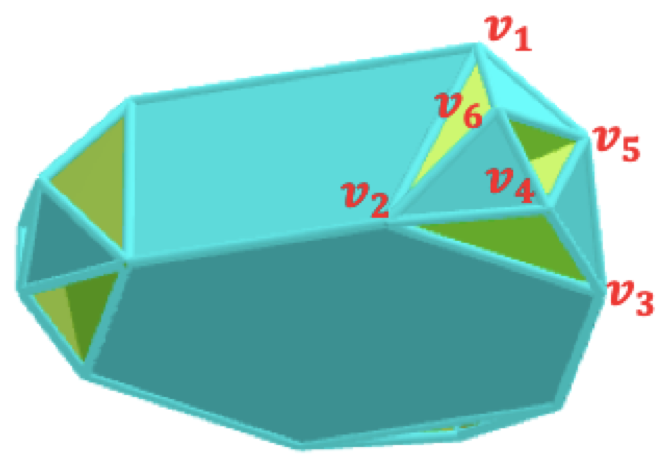}\hfill
    \includegraphics[width=0.2\linewidth]{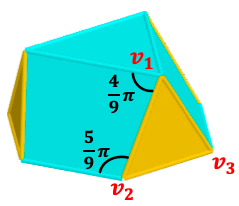} \hfill
    \includegraphics[width=0.4\linewidth]{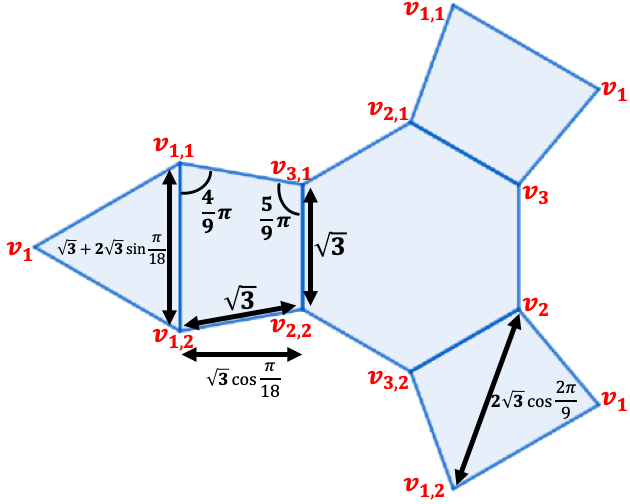} \hfill
    \includegraphics[width=0.3\linewidth]{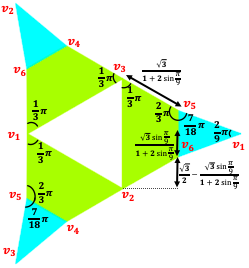}
    \caption{Top row, three panels from left: a non-orientable genus-$3$ CCP $Q^3_{18}$. 
    Top right: a polyhedron $S$ with the topology of a sphere. 
    Bottom row: nets of $S$ (without the lateral triangular faces) and the projective plane $R(r,h)$.}
    \label{fig:n_g3_v18}
\end{figure}

For $g = 3$, a CCP can be constructed by forming the connected sum of a trapezoidal polyhedron $S$, topologically equivalent to a sphere, with three projective planes $R(r,h)$. 
The polyhedron $S$ consists of a regular hexagonal base and a regular triangular top, connected by three regular triangles and isosceles trapezoids, as depicted in \Cref{fig:n_g3_v18} (Top right).
The coordinates of its base vertices are given by
\[
v_1 = \left( \frac{3}{2} - \sqrt{\frac{9}{4} - h_2^2},\, 0,\, h_2 \right), \quad
v_2 = \left( \frac{3}{2},\, -\frac{\sqrt{3}}{2},\, 0 \right), \quad
v_3 = \left( \frac{3}{2},\, \frac{\sqrt{3}}{2},\, 0 \right),
\]
where
\[
h_2 = \sqrt{ -4\sin^2{\frac{\pi}{18}} + 2\sin{\frac{\pi}{18}} + 2 }.
\]
The coordinates of the other vertices $v_{1,k}$ are obtained by rotation by $\frac{2k\pi}{3}$ about the $z$-axis.
With this choice of $h_2$, we have
$\angle v_{1}v_{2}v_{3,2}=\frac{5\pi}{9}$.

Let $r = \dfrac{2\sin{\frac{\pi}{9}}}{1 + 2\sin{\frac{\pi}{9}}}$ and $h = \dfrac{ \sqrt{ -4\sin^2{\frac{\pi}{9}} - 2\sin{\frac{\pi}{9}} + 2 } }{ 1 + 2\sin{\frac{\pi}{9}} }$.
The polyhedron $Q^3_{18} = S \# 3R(r, h)$, shown in \Cref{fig:n_g3_v18} (Right), is a non-orientable genus-$3$ CCP. 
Here, the three connected sums are performed along the face $v_1 v_2 v_3$ and its images under rotation about the $z$-axis.
It is straightforward to check that $|V| = 18$, $|E| = 42$, and $|F| = 23$, and hence
\[
2-g = \chi(Q^3_{18}) = |V| - |E| + |F| = -1.
\]
Thus, 
by \eqref{eq:Descartes}, the angular defect at each vertex should be
\[
\delta_v = \frac{2\pi (2 - g)}{18} = -\frac{\pi}{9}.
\]
To verify this, we compute
\begin{align*}
\sqrt{\frac{9}{4} - {h_2}^2} &= \sqrt{ \frac{9}{4} - \left(-4\sin^2{\frac{\pi}{18}} + 2\sin{\frac{\pi}{18}} + 2\right) }\\
& = \sqrt{\left(-2\sin{\frac{\pi}{18}} + \frac{1}{2}\right)^2} = -2\sin{\frac{\pi}{18}} + \frac{1}{2},\\[1em]
\sqrt{r^2-r+1+h^2}
&= \sqrt{\left(\frac{2\sin{\frac{\pi}{9}}}{1+2\sin{\frac{\pi}{9}}}\right)^2 - \frac{2\sin{\frac{\pi}{9}}}{1+2\sin{\frac{\pi}{9}}} + 1 + \left(\frac{\sqrt{-4\sin^2{\frac{\pi}{9}}-2\sin{\frac{\pi}{9}}+2}}{1+2\sin{\frac{\pi}{9}}}\right)^2}\\
& = \frac{\sqrt{\left(2\sin{\frac{\pi}{9}}\right)^2-\left(2\sin{\frac{\pi}{9}}\right)\left(1+2\sin{\frac{\pi}{9}}\right)+\left(1+2\sin{\frac{\pi}{9}}\right)^2+\left(-4\sin^2{\frac{\pi}{9}}-2\sin{\frac{\pi}{9}}+2\right)}}{1+2\sin{\frac{\pi}{9}}}\\
& = \frac{\sqrt{3}}{1+2\sin{\frac{\pi}{9}}}\\
\end{align*}

% \begin{align*}
% | v_{1,1} v_{1,2} | &= \Bigg|
% \left(
% -\frac{1}{2}\left(\frac{3}{2} - \sqrt{\frac{9}{4} - {h_2}^2}\right),\,
% \frac{\sqrt{3}}{2}\left(\frac{3}{2} - \sqrt{\frac{9}{4} - {h_2}^2}\right),\,
% h_2
% \right), \\
% &\hspace{3em}
% \left(
% -\frac{1}{2}\left(\frac{3}{2} - \sqrt{\frac{9}{4} - {h_2}^2}\right),\,
% -\frac{\sqrt{3}}{2}\left(\frac{3}{2} - \sqrt{\frac{9}{4} - {h_2}^2}\right),\,
% h_2
% \right)
% \Bigg| \\
% &= \frac{3\sqrt{3}}{2} - \sqrt{3}\sqrt{ \frac{9}{4} - {h_2}^2 }\\
% &= \frac{3\sqrt{3}}{2} - \sqrt{3}\left(-2\sin{\frac{\pi}{18}} + \frac{1}{2}\right)\\
% &= \sqrt{3}\left(1+2\sin{\frac{\pi}{18}}\right), \\[1em]
% | v_{3,1} v_{2,2} | &= \left| \left(-\frac{3}{2},\, \frac{\sqrt{3}}{2},\, 0\right),\; \left(-\frac{3}{2},\, -\frac{\sqrt{3}}{2},\, 0\right) \right|=\sqrt{3}, \\[1em]
% | v_{1,2} v_{2,2} | &= \Bigg|
% \left(-\frac{1}{2}\left(\frac{3}{2} - \sqrt{\frac{9}{4} - {h_2}^2}\right),\, -\frac{\sqrt{3}}{2}\left(\frac{3}{2} - \sqrt{\frac{9}{4} - {h_2}^2}\right),\, h_2\right), \\
% &\hspace{3em}
% \left(-\frac{3}{2},\, -\frac{\sqrt{3}}{2},\, 0\right)
% \Bigg| = \sqrt{3}, \\[1em]
% | v_{1,1} v_{2,2} | &= \Bigg|
% \left(
% -\frac{1}{2}\left(\frac{3}{2} - \sqrt{\frac{9}{4} - {h_2}^2}\right),\,
% \frac{\sqrt{3}}{2}\left(\frac{3}{2} - \sqrt{\frac{9}{4} - {h_2}^2}\right),\,
% h_2
% \right), \\
% &\hspace{3em}
% \left(-\frac{3}{2},\, -\frac{\sqrt{3}}{2},\, 0\right)
% \Bigg| = \sqrt{\frac{15}{2}-3\sqrt{\frac{9}{4}-{h_2}^2}}\\
% &= \sqrt{\frac{15}{2}-3\left(-2\sin{\frac{\pi}{18}} + \frac{1}{2}\right)}\\
% &= \sqrt{6\left(\sin{\frac{\pi}{18}}+1\right)}.
% \end{align*}
Putting  $k = \frac{3}{2} - \sqrt{\frac{9}{4} - {h_2}^2}$, 
\begin{align*}
| v_{1,1} v_{1,2} | &= \left| \left(-\frac{k}{2},\, \frac{k\sqrt{3}}{2},\, h_2\right),\, \left(-\frac{k}{2},\, -\frac{k\sqrt{3}}{2},\, h_2\right) \right| \\
&= k\sqrt{3} = \sqrt{3}\left(1+2\sin\frac{\pi}{18}\right), \\[1em]
| v_{3,1} v_{2,2} | &= \left| \left(-\frac{3}{2},\, \frac{\sqrt{3}}{2},\, 0\right),\, \left(-\frac{3}{2},\, -\frac{\sqrt{3}}{2},\, 0\right) \right| = \sqrt{3}, \\[1em]
| v_{1,2} v_{2,2} | &= \left| \left(-\frac{k}{2},\, -\frac{k\sqrt{3}}{2},\, h_2\right),\, \left(-\frac{3}{2},\, -\frac{\sqrt{3}}{2},\, 0\right) \right| \\
&= \sqrt{3 - 3k + k^2 + {h_2}^2} = \sqrt{3}, \\[1em]
| v_{1,1} v_{2,2} | &= \left| \left(-\frac{k}{2},\, \frac{k\sqrt{3}}{2},\, h_2\right),\, \left(-\frac{3}{2},\, -\frac{\sqrt{3}}{2},\, 0\right) \right| \\
&= \sqrt{3(k+1)} = \sqrt{6\left(\sin\frac{\pi}{18}+1\right)}.
\end{align*}

Now, the relevant angles are given by
\begin{align*}
\angle v_1 v_5 v_6 &= \arccos{\frac{\sqrt{3}r}{2\sqrt{r^2 - r + 1 + h^2}}}
    = \arccos{\frac{\sqrt{3}\frac{2\sin{\frac{\pi}{9}}}{1+2\sin{\frac{\pi}{9}}}}{2\frac{\sqrt{3}}{1+2\sin{\frac{\pi}{9}}}}}\\
    &= \arccos{\left(\sin{\frac{\pi}{9}}\right)}
    = \arccos{\left(\cos{\frac{7\pi}{18}}\right)}
    = \frac{7\pi}{18}\\
\angle v_5v_1v_6 &= \pi - (\angle v_1 v_5 v_6 + \angle v_1 v_6 v_5)
    = \pi - (\angle v_1 v_5 v_6 + \angle v_1 v_5 v_6)\\
    &= \pi - 2\,\angle v_1 v_5 v_6
    = \pi - 2\,\frac{7\pi}{18}
    = \frac{2\pi}{9}\\
\angle v_5 v_3 v_2 &= \arccos{\frac{\sqrt{3} \left(1 - r\right)}{2\sqrt{r^2 - r + 1 + h^2}}}
    = \arccos{\frac{\sqrt{3} \left(1 - \frac{2\sin{\frac{\pi}{9}}}{1+2\sin{\frac{\pi}{9}}}\right)}{2\frac{\sqrt{3}}{1+2\sin{\frac{\pi}{9}}}}}\\
    &= \arccos{\frac{1+2\sin{\frac{\pi}{9}}-2\sin{\frac{\pi}{9}}}{2}}
    = \frac{\pi}{3}\\
\angle v_3 v_5 v_6 &= \pi - \angle v_5 v_3 v_2
    = \pi - \frac{\pi}{3}
    = \frac{2\pi}{3}.
\end{align*}

\begin{align*}
\angle v_{2,2} v_{1,2} v_{1,1} &= \arccos{ \frac{ |v_{1,1}v_{1,2}|^2 + |v_{1,2}v_{2,2}|^2 - |v_{1,1}v_{2,2}|^2}{ 2|v_{1,1}v_{1,2}||v_{1,2}v_{2,2}| } }\\
    &= \arccos{ \frac{ \left(\sqrt{3}\left(1+2\sin{\frac{\pi}{18}}\right)\right)^2 + 3 - 6\left(\sin{\frac{\pi}{18}}+1\right)}{ 2\sqrt{3}\left(1+2\sin{\frac{\pi}{18}}\right)\sqrt{3} } }\\
    &= \arccos{ \frac{6\left(1+2\sin{\frac{\pi}{18}}\right)\sin{\frac{\pi}{18}}}{ 6\left(1+2\sin{\frac{\pi}{18}}\right) } }\\
    &= \arccos{\left(\sin{\frac{\pi}{18}}\right)}\\
    &= \arccos{\left(\cos{\left(\frac{\pi}{2}-\frac{\pi}{18}\right)}\right)}\\
    &= \frac{4\pi}{9},\\[1em]
\angle v_{1,2} v_{2,2} v_{3,1} &= \pi - \angle v_{2,2} v_{1,2} v_{1,1}
    = \pi - \frac{4\pi}{9}
    = \frac{5\pi}{9}
\end{align*}

and the angular defect at the main vertices can be computed as
\begin{align*}
\delta_{v_1} &= 2\pi - \big(
\angle v_{2} v_{1} v_{1,2} +
\angle v_{1,2} v_{1} v_{1,1} +
\angle v_{1,1} v_{1} v_{3} \\
&\hspace{2.6em} +\; \angle v_3 v_1 v_6 +
\angle v_6 v_1 v_5 +
\angle v_5 v_1 v_2
\big) \\
    &= 2\pi - \big(
    \angle v_{2,2} v_{1,2} v_{1,1} +
    \angle v_{1,2} v_{1} v_{1,1} +
    \angle v_{2,2} v_{1,2} v_{1,1} \\
    &\hspace{2.6em} +\; \angle v_5 v_3 v_2 +
    \angle v_5 v_1 v_6 +
    \angle v_5 v_3 v_2
    \big) \\
    &= 2\pi - \left(
    \frac{4\pi}{9} + \frac{\pi}{3} + \frac{4\pi}{9}
    + \frac{\pi}{3} + \frac{2\pi}{9} + \frac{\pi}{3}
    \right) = -\frac{\pi}{9}, \\[1em]
\delta_{v_2} &= 2\pi - \big(
\angle v_{1} v_{2} v_{3,2} +
\angle v_{3,2} v_{2} v_{3} +
\angle v_3 v_2 v_6 \\
&\hspace{2.6em} +
\angle v_6 v_2 v_4 +
\angle v_4 v_2 v_1
\big) \\
    &= 2\pi - \big(
    \angle v_{1,2} v_{2,2} v_{3,1} +
    \angle v_{3,2} v_{2} v_{3} +
    \angle v_5 v_3 v_2 \\
    &\hspace{2.6em} +
    \angle v_5 v_1 v_6 +
    \angle v_5 v_3 v_2
    \big) \\
    &= 2\pi - \left(
    \frac{5\pi}{9} + \frac{2\pi}{3} + \frac{\pi}{3}
    + \frac{2\pi}{9} + \frac{\pi}{3}
    \right) = -\frac{\pi}{9}, \\[1em]
\delta_{v_6} &= 2\pi - \big(
\angle v_1 v_6 v_4 +
\angle v_4 v_6 v_2 +
\angle v_2 v_6 v_5 +
\angle v_5 v_6 v_1
\big) \\
    &= 2\pi - \big(
    \angle v_3 v_5 v_6 +
    \angle v_1 v_5 v_6 +
    \angle v_3 v_5 v_6 +
    \angle v_1 v_5 v_6
    \big) \\
    &= 2\pi - \left(
    \frac{2\pi}{3} + \frac{7\pi}{18} + \frac{2\pi}{3} + \frac{7\pi}{18}
    \right) = -\frac{\pi}{9}.
\end{align*}

\begin{figure}[htbp]
    \centering
    \includegraphics[width=0.2\linewidth]{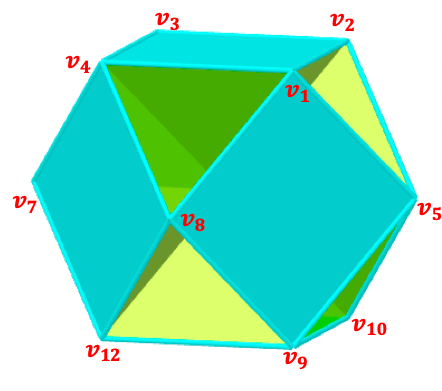}
    \includegraphics[width=0.2\linewidth]{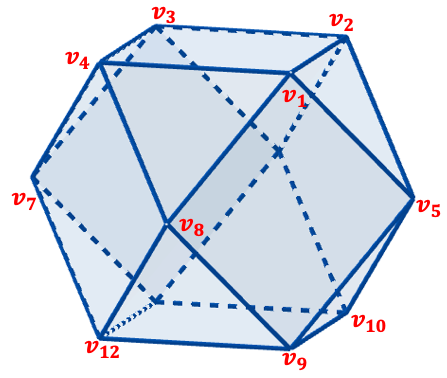}
    \includegraphics[width=0.2\linewidth]{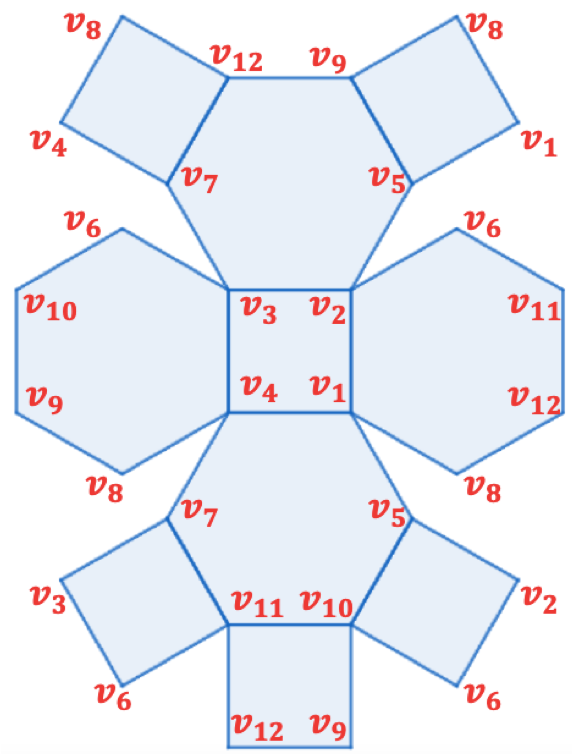}
    \caption{The cubohemioctahedron $Q^4_{12}$, a non-orientable genus-$4$ CCP.}
    \label{fig:n_g4_v12}
\end{figure}

For genus $g=4$, the cubohemioctahedron $Q^4_{12}$, shown in \Cref{fig:n_g4_v12}, is a CCP.
In this polyhedron, two regular hexagons and two squares meet at each vertex, yielding an angular defect of $-\frac{\pi}{3}$ at every vertex.
It is easy to see that $|V| = 12$, $|E| = 24$, and $|F| = 10$, and thus
\[
2-g = \chi(Q^4_{12}) = |V| - |E| + |F| = -2.
\]

For genus $2n + 3$ with $n > 0$, we apply the drilling operation to $Q^3_{18}$ along a segment $\ell$ connecting the top and bottom faces in \Cref{fig:n_g3_v18}, obtaining
\[
Q^{2n+3}_{18 + 36n} = Q^3_{18} \# n\Pi(\ell, 18),
\]
where $18=-\frac{|V|}{\chi}$ (by \eqref{eq:choice-of-n}).

For genus $2n + 4$ with $n > 0$, we apply the drilling operation to $Q^4_{12}$ along a segment $\ell$ connecting the top and bottom faces in \Cref{fig:n_g4_v12}, obtaining
\[
Q^{2n+4}_{12 + 12n} = Q^4_{12} \# n\Pi(\ell, 6),
\]
where $6=-\frac{|V|}{\chi}$ (by \eqref{eq:choice-of-n}).

%%%%%%%%%%%%%%%%%%%%%
\section{Constant Curvature Polyhedra with Fewer Vertices}
In the study of polyhedra, constructions with minimal elements have been actively explored, as exemplified by the polyhedra of Császár and Szilassi \cite{ringel1955man,minimal}. The Szilassi polyhedron, with 14 vertices, 21 edges, and 7 faces, represents the minimum number of vertices for a toroidal polyhedron. Its dual, the Császár polyhedron, has 7 vertices, 21 edges, and 14 faces, and achieves the minimum number of faces for a toroidal polyhedron.

In this section, we consider CCPs with a small number of vertices. 
Our approach is to provide explicit constructions of CCPs, thereby establishing upper bounds for the minimal vertex count in each case. Determining lower bounds appears to be significantly more challenging; for instance, in the case of flat torus topology, no theoretical lower bound results appear to be available in the literature. The construction shown in \Cref{fig:o_g1} has 9 vertices and appears to be minimal among CCPs with flat torus topology, in contrast to the Császár polyhedron which has 7 vertices for torus topology but without the constant angular defect constraint that defines CCPs.
\Cref{tab:vertex} summarises the currently known CCPs with a small number of vertices.
Most of them are new and constructed in this section and appendix.

\begin{table}[htbp]
    \centering
\caption{Known CCPs by genus and orientability}
\label{tab:vertex}
\begin{tabular}{@{}lllc@{}}
\toprule
\textbf{Orientability} & \textbf{Self-intersection} & \textbf{Genus ($g$)} & \textbf{Vertices} \\
\midrule
\multirow{8}{*}{\rotatebox{90}{Orientable}} 
& \multirow{6}{*}{\rotatebox{0}{Without}} 
& $g = 0$ & 4(regular tetrahedron) \\
& & $g = 1$ & 9(\Cref{fig:o_g1}) \\
& & $g = 2, 3$ & $8g$(\Cref{fig:o_g3_v24}) \\
& & $g = 4, 5, 6$ & $7g - 7$(\Cref{fig:o_g5_v28}) \\
& & $g \geq 7$ & $6g$(\Cref{fig:o_g8_v48}) \\
\cmidrule{2-4}
& \multirow{1}{*}{\rotatebox{0}{Allowed}} 
 & $g \geq 0$ & $2g+4$ \\
\midrule
\multirow{8}{*}{\rotatebox{90}{Non-orientable}} 
& \multirow{8}{*}{\rotatebox{0}{Allowed}} 
& $g = 1$ & 6(\Cref{fig:n_g1_v6}) \\
& & $g = 3, 5, 7, 9, 11$ & $5g$(\Cref{fig:n_g5_v25}) \\
& & $g \geq 13$ and odd & $7g - 14$(\Cref{fig:n_g5_v25}) \\
\cmidrule{3-4}
& & $g = 2$ & 9(\Cref{fig:n_g2_v9}) \\
& & $g = 4, 6$ & $6g - 12$(\Cref{fig:n_g4_v12}) \\
& & $g = 14$ & 30(\Cref{fig:n_g14_v30}) \\
& & other even numbers & $4g - 8$(\Cref{fig:n_g8_v24}) \\
\bottomrule
\end{tabular}
\end{table}

\begin{remark}
    In our definition of a polyhedron (see \Cref{dfn:polyhedron}), we exclude those having stellar polygons as faces. There exist stellar polyhedra with a small number of vertices that realise non-orientable surfaces; for example,
    Ditrigonal Dodecadodecahedron (orientable, $g=9$, $|V|=20$), 
    Small dodecahemicosahedron (non-orientable, $g=10$, $|V|=30$),
    Great Dodecahemioctahedron (non-orientable, $g=14$, $|V|=30$),
    Great rhombidodecahedron (non-orientable, $g=20$, $|V|=60$),
    and Great dodecicosahedron (non-orientable, $g=30$, $|V|=60$).
\end{remark}

For orientable CCP without self-intersections and for non-orientable CCP, the known examples listed in \Cref{tab:vertex} are likely not minimal in terms of vertex count. 
But we conjecture:
\begin{conjecture}
\ 
    \begin{enumerate}[label={[\roman*]}]
        \item For a CCP of a topological torus without self intersection, the minimal number of vertices is $9$.
        \item For orientable CCP of genus~$g$, if self-intersections are permitted, then the minimal number of vertices is $2g + 4$.
    \end{enumerate}
\end{conjecture}
In what follows, we shall construct orientable CCPs of genus~$g$ with self-intersections having $2g+4$ vertices, providing supportive evidence for conjecture (ii).

\subsection{Orientable, with self-intersection, $2g+4$ vertices}
\begin{theorem}\label{thm:minimal}
    For every integer $g \ge 1$ there exists an orientable genus-$g$ CCP with $2g+4$ vertices and self-intersection.
\end{theorem}

\begin{figure}[htbp]
    \centering
    \includegraphics[width=0.25\linewidth]{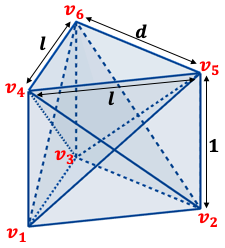}\hfill
    \includegraphics[width=0.25\linewidth]{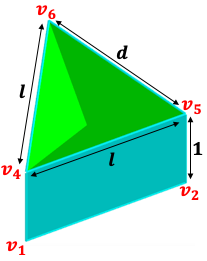}\hfill
    \includegraphics[width=0.35\linewidth]{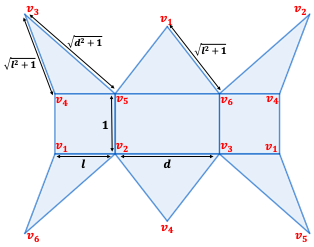}
    \caption{The genus-$1$ polyhedron $T(l,d)$.}
    \label{fig:o_g1_v6}
\end{figure}

We begin with the genus–$1$ polyhedron $T(l,d)$ in \Cref{fig:o_g1_v6}, with vertices $v_1,\dots,v_6$ and parameters $l,d>0$ controlling edge lengths.  
Setting $d=l$ yields a CCP with zero angle defect at every vertex.  
The polyhedron $T(l,d)$ will serve as a building block for higher–genus CCPs.

\begin{figure}[htbp]
    \centering
    \begin{minipage}[b]{0.45\linewidth}
        \includegraphics[width=\linewidth]{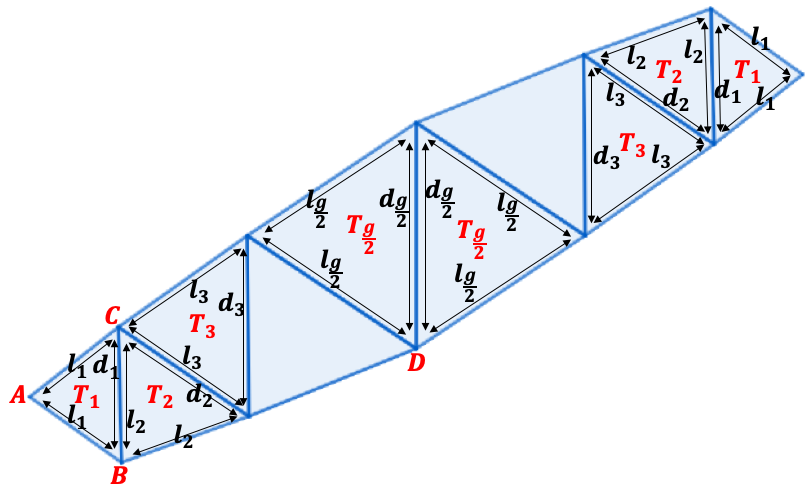}
        \subcaption*{even $g$}
    \end{minipage}
    \begin{minipage}[b]{0.45\linewidth}
        \includegraphics[width=\linewidth]{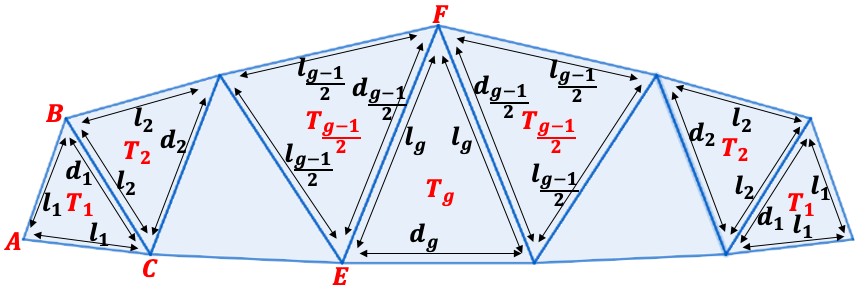}
        \subcaption*{odd $g$}
    \end{minipage}
    \caption{Top view of the gluing pattern for an orientable CCP with $2g+4$ vertices. Subscripts label the genus–1 tori $T_k$.}
    \label{fig:o_v2g+4_top}
\end{figure}

For $g\ge2$ we glue $g$ copies of $T(l,d)$ (with possibly different parameter pairs $(l,d)$) along rectangular side faces, as indicated in \Cref{fig:o_v2g+4_top}, obtaining a closed orientable polyhedral surface.  
The configuration admits a reflection in a plane which interchanges the triangles $v_1v_2v_3$ and $v_4v_5v_6$, so in what follows we restrict attention to the ``upper'' side containing the triangle $v_4v_5v_6$.

In the glued polyhedron we distinguish the following types of vertices (see \Cref{fig:o_v2g+4_top}):

\begin{description}[labelsep=0.6em,leftmargin=1.9em]
  \item[($A$)] A vertex contained in a single block, namely $v_4$ of $T_1$.
  \item[($B$)] A vertex at which two blocks meet; for example, $v_5$ (or $v_6$) of $T_1$, identified with $v_4$ of $T_{2}$.
  \item[($C$)] A vertex at which three blocks meet. For instance, the blocks $T_{k-1}$, $T_k$, $T_{k+1}$ contribute $v_6$, $v_6$, $v_4$ (or $v_5$, $v_5$, $v_4$), which are identified to a single vertex.
  \item[($D$)] (Appears only when $g$ is even.) A vertex where the two copies of $T_{g/2}$ meet; if $g\ge4$ the block $T_{g/2-1}$ also contributes.
  \item[($E$)] (Appears only when $g$ is odd.) A vertex where $T_{(g-1)/2}$ and $T_{(g-1)/2-1}$ meet; if $g\ge5$ the block $T_{(g-1)/2-2}$ also contributes.
  \item[($F$)] (Appears only when $g$ is odd.) A vertex where three blocks meet: $v_4$ of the central $T_g$ is identified with $v_6$ (or $v_5$) coming from each of the two copies of $T_{(g-1)/2}$ attached to it.
\end{description}

We shall verify that the angle defect $\delta$ takes the same value at all vertices of types $A,B,C,F$
by choosing appropriate parameters.
The remaining end types $D,E$ will then be forced to share this common value by Descartes' theorem.

\medskip

We introduce the coefficients
\[
a_{k,g}:=\frac{4\Bigl(3k+1-\frac{1}{(-2)^k}\Bigr)}{3(g+2)}\pi
\qquad (k\in\mathbb{Z}_{\ge1}),
\]
which are strictly increasing in $k$ and satisfy $0<a_{k,g}<2\pi$.  
Denote by $v_{i|T_k}$ the vertex $v_i$ on the block $T_k$.  
Suppose that there exist parameters $(l_k,d_k)$ so that for each block $T_k$ with $k=1,\dots,\lfloor g/2\rfloor$ the sum of the three angles at $v_4$ satisfies
\[
\angle v_{6|T_k}v_{4|T_k}v_{2|T_k}
+\angle v_{2|T_k}v_{4|T_k}v_{3|T_k}
+\angle v_{3|T_k}v_{4|T_k}v_{5|T_k}
=3\pi-a_{k,g},
\]
and, when $g$ is odd, the terminal block $T_g$ additionally satisfies
\[
\angle v_{6|T_g}v_{4|T_g}v_{2|T_g}
+\angle v_{2|T_g}v_{4|T_g}v_{3|T_g}
+\angle v_{3|T_g}v_{4|T_g}v_{5|T_g}
=3\pi-a_{\frac{g-1}{2},g}-\frac{6}{g+2}\pi.
\]

Then,
\begin{align*}
&\angle v_{4|T_k}v_{5|T_k}v_{3|T_k}
 +\angle v_{3|T_k}v_{5|T_k}v_{1|T_k}
 +\angle v_{1|T_k}v_{5|T_k}v_{6|T_k}\\
&\quad=
 \angle v_{4|T_k}v_{5|T_k}v_{3|T_k}
 +\angle v_{5|T_k}v_{3|T_k}v_{4|T_k}
 +\angle v_{4|T_k}v_{2|T_k}v_{3|T_k}
 \qquad\text{(by the local configuration in \Cref{fig:o_g1_v6})}\\
&\quad=
 \bigl(\angle v_{4|T_k}v_{5|T_k}v_{3|T_k}
      +\angle v_{5|T_k}v_{3|T_k}v_{4|T_k}\bigr)
 +\frac12\bigl(\angle v_{4|T_k}v_{2|T_k}v_{3|T_k}
               +\angle v_{4|T_k}v_{3|T_k}v_{2|T_k}\bigr)\\
&\quad=
 \bigl(\pi-\angle v_{3|T_k}v_{4|T_k}v_{5|T_k}\bigr)
 +\frac12\bigl(\pi-\angle v_{2|T_k}v_{4|T_k}v_{3|T_k}\bigr)
 \qquad\text{(angle sums in triangles $v_3v_4v_5$ and $v_2v_3v_4$)}\\
&\quad=
 \frac32\pi-\angle v_{3|T_k}v_{4|T_k}v_{5|T_k}
             -\frac12\angle v_{2|T_k}v_{4|T_k}v_{3|T_k}\\
&\quad=
 \frac32\pi-\frac12\Bigl(
   \angle v_{6|T_k}v_{4|T_k}v_{2|T_k}
  +\angle v_{2|T_k}v_{4|T_k}v_{3|T_k}
  +\angle v_{3|T_k}v_{4|T_k}v_{5|T_k}\Bigr)
 \qquad\text{(by $v_5\leftrightarrow v_6$ symmetry)}\\
&\quad=
 \frac12\Bigl(3\pi-
   \bigl[\angle v_{6|T_k}v_{4|T_k}v_{2|T_k}
        +\angle v_{2|T_k}v_{4|T_k}v_{3|T_k}
        +\angle v_{3|T_k}v_{4|T_k}v_{5|T_k}\bigr]\Bigr)\\
&\quad=
 \frac12\,a_{k,g},
\end{align*}
where in the last line we used the imposed relation
\[
\angle v_{6|T_k}v_{4|T_k}v_{2|T_k}
+\angle v_{2|T_k}v_{4|T_k}v_{3|T_k}
+\angle v_{3|T_k}v_{4|T_k}v_{5|T_k}
=3\pi-a_{k,g}.
\]
By the same argument with $v_5$ replaced by $v_6$ we also obtain
\[
\angle v_{4|T_k}v_{6|T_k}v_{2|T_k}
+\angle v_{2|T_k}v_{6|T_k}v_{1|T_k}
+\angle v_{1|T_k}v_{6|T_k}v_{5|T_k}
=\frac12\,a_{k,g}.
\]

When $g$ is odd and $T_g$ is the terminal block, the same computations yield the corresponding identities with right–hand side
\[
\frac12\!\left(a_{\frac{g-1}{2},g}+\frac{6}{g+2}\pi\right).
\]
These local identities are the only ingredients needed in what follows to sum the angles at vertices of types $A,B,C,F$ and to verify that they all share the same defect $\delta$.  
%The construction of parameters $(l_k,d_k)$ ensuring the prescribed $v_4$–sum relations will be given afterwards.

\medskip

We first compute the defect at a type $A$ vertex.  
Around the $v_4$–corner of $T_1$ there are two right angles (coming from the two rectangles) and the three–angle sum at $v_4$, so
\[
\delta_A
=2\pi-\Bigl(\frac12\pi+\bigl(3\pi-a_{1,g}\bigr)+\frac12\pi\Bigr)
=-2\pi+a_{1,g}
=-2\pi+\frac{6}{g+2}\pi .
\]

Next, consider a type $B$ vertex where two consecutive blocks $T_1$ and $T_2$ meet across a glued rectangle.  
Here we collect two right angles, the ``$v_5$–side'' triple–angle in $T_1$, and the $v_4$–sum in $T_2$:
\[
\delta_B
=2\pi-\Bigl(\frac12\pi+\frac12\,a_{1,g}+\bigl(3\pi-a_{2,g}\bigr)+\frac12\pi\Bigr)
=-2\pi-\frac12 a_{1,g}+a_{2,g}
=-2\pi+\frac{6}{g+2}\pi .
\]

For a type $C$ vertex, three consecutive blocks $T_{k-1},T_k,T_{k+1}$ meet.  
Summing the two right angles together with the ``$v_6$–side'' triple–angle in $T_{k-1}$, the ``$v_5$–side'' triple–angle in $T_k$, and the $v_4$–sum in $T_{k+1}$ gives
\begingroup\setlength{\jot}{.25ex}
\begin{align*}
\delta_C
&=2\pi-\Bigl(\frac12\pi+\frac12 a_{k-1,g}+\frac12 a_{k,g}
+\bigl(3\pi-a_{k+1,g}\bigr)+\frac12\pi\Bigr)\\
&=-2\pi-\frac12 a_{k-1,g}-\frac12 a_{k,g}+a_{k+1,g}
=-2\pi+\frac{6}{g+2}\pi.
\end{align*}
\endgroup

Finally, assume $g$ is odd and consider a type $F$ vertex, where the central block $T_g$ meets the two copies of $T_{(g-1)/2}$.  
Using two right angles, the ``$v_6$– or $v_5$–side'' triple–angle in each $T_{(g-1)/2}$, and the $v_4$–sum in $T_g$, we obtain
\begingroup\setlength{\jot}{.25ex}
\begin{align*}
\delta_F
&=2\pi-\Bigl(\frac12\pi+\frac12 a_{\frac{g-1}{2},g}
+\bigl(3\pi-a_{\frac{g-1}{2},g}-\tfrac{6}{g+2}\pi\bigr)
+\frac12 a_{\frac{g-1}{2},g}+\frac12\pi\Bigr)\\
&=-2\pi+\frac{6}{g+2}\pi .
\end{align*}
\endgroup
Thus every vertex of types $A,B,C,F$ has the same defect
\[
\delta=-\frac{2g-2}{g+2}\,\pi=-2\pi+\frac{6}{g+2}\pi .
\]

\medskip

We now consider the end types $D$ and $E$.  
Starting from one torus with $(V,E,F)=(6,15,9)$, each connected sum along a rectangle (as in \Cref{fig:o_v2g+4_top}) changes the counts by
\[
\Delta V=+2,\qquad \Delta E=+11,\qquad \Delta F=+7,
\]
so after assembling a genus–$g$ surface we obtain
\[
V=2g+4,\qquad E=11g+4,\qquad F=7g+2,
\]
and hence $V-E+F=2-2g$.  
By Descartes' formula,
\[
\frac{1}{V}\sum_{v}\delta_v=\frac{2\pi\chi}{V}
=\frac{2\pi(2-g)}{2g+4}
=-\frac{2g-2}{g+2}\,\pi ,
\]
which coincides with the common value computed above for types $A,B,C,F$.  
Therefore the remaining vertex types $D,E$ must also have defect $\delta$ (otherwise the total sum of defects would not match), and the glued surface is indeed a CCP.

\medskip

It remains to construct parameters $(l_k,d_k)$ realising the imposed $v_4$–sum relations.

Define, for $l>0$,
\[
f_l(d):=\angle v_{6}v_{4}v_{2}
+\angle v_{2}v_{4}v_{3}
+\angle v_{3}v_{4}v_{5}.
\]
From the explicit coordinates in \Cref{fig:o_g1_v6} we obtain
\[
\angle v_{5}v_{4}v_{3}
=\arccos\!\frac{2l^{2}-d^{2}}{2l\sqrt{l^{2}+1}},\qquad
\angle v_{3}v_{4}v_{2}
=\arccos\!\frac{2(l^{2}+1)-d^{2}}{2(l^{2}+1)},\qquad
\angle v_{2}v_{4}v_{6}
=\angle v_{5}v_{4}v_{3},
\]
and hence
\[
f_l(d)
=2\arccos\!\frac{2l^{2}-d^{2}}{2l\sqrt{l^{2}+1}}
+\arccos\!\frac{2(l^{2}+1)-d^{2}}{2(l^{2}+1)}.
\]
In particular,
\[
f_l(0)=\pi-2\arctan l,\qquad
f_l(l)=\pi,\qquad
f_l(2l)=\pi+4\arctan l,
\]
and 
\[
\lim_{l\to\infty}f_l(0)=0,\qquad
\lim_{l\to\infty}f_l(2l)=3\pi .
\]

Choose $l_1>0$ sufficiently large so that
\[
f_{l_1}(2l_1)>3\pi-a_{1,g},
\qquad
\text{and, if $g$ is odd,}\quad
f_{l_1}(0)<3\pi-a_{\frac{g-1}{2},g}-\frac{6}{g+2}\pi .
\]
Since $f_{l_1}(l_1)=\pi<3\pi-a_{1,g}<f_{l_1}(2l_1)$ and $f_{l_1}$ is continuous and strictly increasing on $[l_1,2l_1]$, the intermediate value theorem gives a unique $d_1\in(l_1,2l_1)$ such that
\[
f_{l_1}(d_1)=3\pi-a_{1,g}.
\]
Inductively define $l_{k+1}:=d_k$.  
Because $a_{k+1,g}>a_{k,g}$ while
\[
f_{l_{k+1}}(l_{k+1})=\pi,\qquad
f_{l_{k+1}}(2l_{k+1})=\pi+4\arctan l_{k+1}>\pi+4\arctan l_k=f_{l_k}(2l_k),
\]
we have
\[
\pi=f_{l_{k+1}}(l_{k+1})<3\pi-a_{k+1,g}<f_{l_{k+1}}(2l_{k+1}),
\]
so there exists a unique $d_{k+1}\in(l_{k+1},2l_{k+1})$ with
\[
f_{l_{k+1}}(d_{k+1})=3\pi-a_{k+1,g}.
\]
This constructs $(l_k,d_k)$ for all $k=1,\dots,\lfloor g/2\rfloor$.  
If $g$ is odd, we finally set $l_g:=d_{\frac{g-1}{2}}$ and use
\[
f_{l_g}(0)<3\pi-a_{\frac{g-1}{2},g}-\frac{6}{g+2}\pi
< f_{l_g}(2l_g),
\]
to obtain a unique $d_g\in(0,2l_g)$ satisfying
\[
f_{l_g}(d_g)=3\pi-a_{\frac{g-1}{2},g}-\frac{6}{g+2}\pi.
\]
Therefore the imposed $v_4$–sum relations hold for all blocks.

\begin{figure}[htbp]
    \centering
    \includegraphics[width=0.5\linewidth]{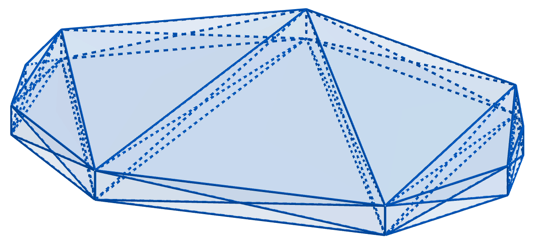}\hfill
    \includegraphics[width=0.5\linewidth]{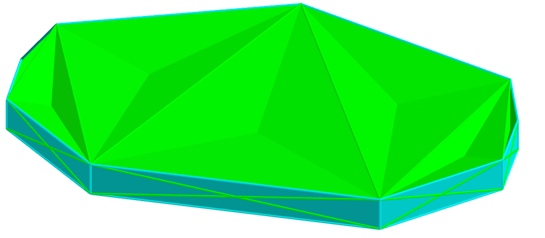}
    \caption{
        An orientable genus-$7$ CCP with 18 vertices and self-intersection.
    }
    \label{fig:o_g7_v18}
\end{figure}

\begin{figure}[htbp]
    \centering
    \includegraphics[width=0.5\linewidth]{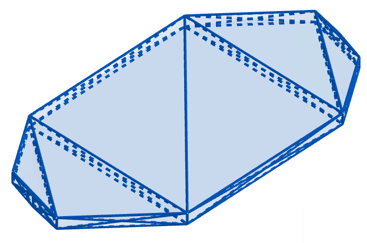}\hfill
    \includegraphics[width=0.5\linewidth]{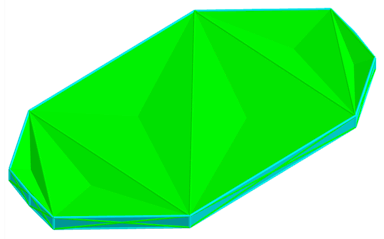}
    \caption{
        An orientable genus-$8$ CCP with 20 vertices and self-intersection.
    }
    \label{fig:o_g8_v20}
\end{figure}

\begin{example}
Figures \Cref{fig:o_g7_v18} and \Cref{fig:o_g8_v20} show the orientable genus–7 CCP and the orientable genus–8 CCP, respectively.  
The corresponding parameter values are approximately:
\[
\text{$g=7$:}\quad
l_1=2,\quad l_2=d_1=3.94799\ldots,\quad
l_3=d_2=6.93234\ldots,\quad
l_7=d_3=9.83752\ldots,\quad
d_7=8.30361\ldots,
\]
\[
\text{$g=8$:}\quad
l_1=2,\quad l_2=d_1=3.99386\ldots,\quad
l_3=d_2=7.21534\ldots,\quad
l_4=d_3=11.01272\ldots,\quad
d_4=13.64880\ldots.
\]
\end{example}

\section*{Acknowledgements}
We are grateful to Professor T. Tsuboi for the stimulating discussions. An early draft of his paper \cite{tsuboi2024} served as a significant inspiration for the present work. 
We are also thankful to Professor J. Ito for providing his idea for CCPs with a small number of vertices.
The second-named author acknowledges support from KAKENHI, Grant-in-Aid for Scientific Research (B) 25K00921 and (S) 25H00399.

%%%%%%%%%%%%%%%%%%%%
\appendix
\section{Reducing the number of vertices}
We explore some more constructions of CCPs with fewer vertex counts,
which are summarised in~\Cref{tab:vertex}.

\subsection{Orientable, without self-intersection, $8g$ vertices $(g\ge 2)$}

\begin{figure}[htbp]
\centering
\includegraphics[width=0.2\linewidth]{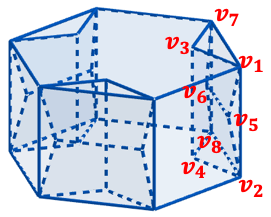}\hfill
\includegraphics[width=0.2\linewidth]{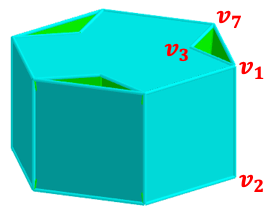}\hfill
\includegraphics[width=0.2\linewidth]{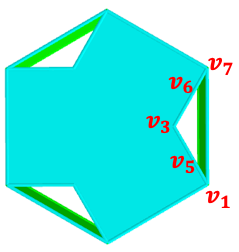}\hfill
\includegraphics[width=0.2\linewidth]{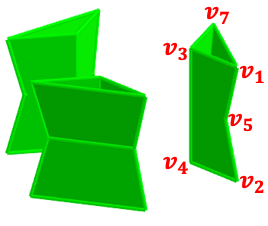}\hfill
\includegraphics[width=0.8\linewidth]{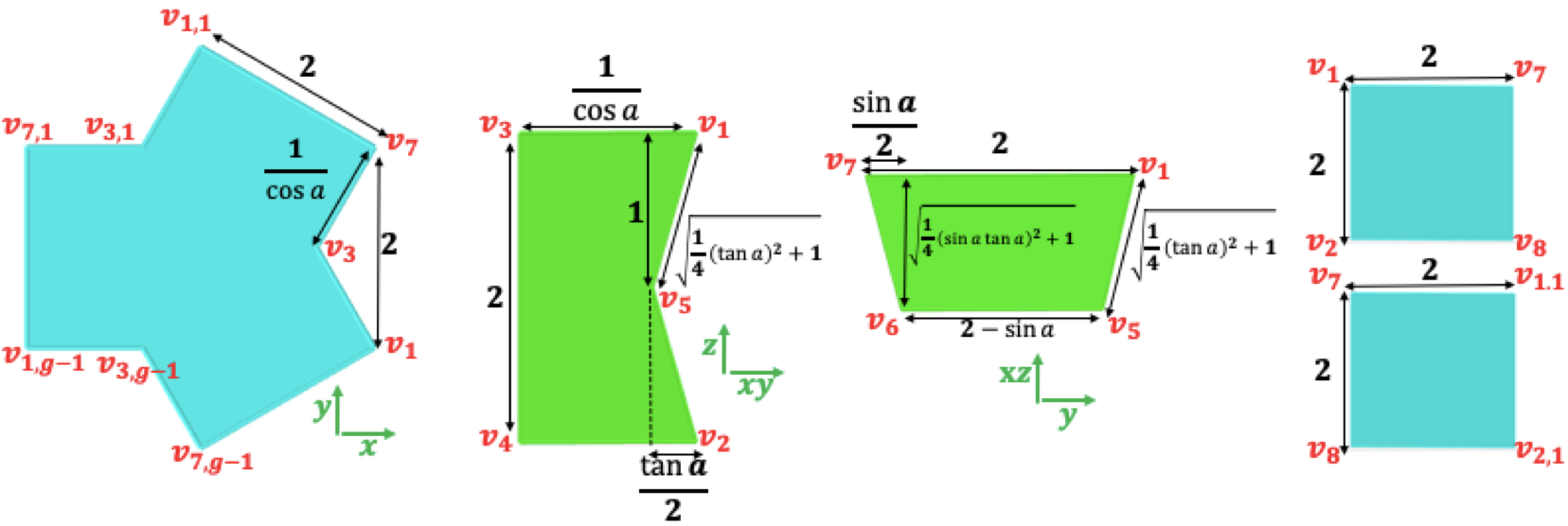}
\caption{
An embedded (without self-intersection)
     orientable CCP of genus $g=3$ with 24 vertices.
}
\label{fig:o_g3_v24}
\end{figure}

Starting from the $g=2$ template in \Cref{fig:o_g2_v24_k}, we reduce the vertex count by
(i) merging the pairs $v_2^{+**}$ with $v_2^{-**}$ and $v_3^{**+}$ with $v_3^{**-}$, and
(ii) placing ``windows'' on every other lateral face of a regular $2g$–gonal prism rather than using the $12$–gonal scheme as in \Cref{fig:o_g3_v24}.
Each handle then contributes only $8$ vertices.

For every $g\ge2$ this yields an embedded orientable CCP with
\[
|V|=8g,\qquad |E|=16g,\qquad |F|=6g+2,\qquad 
\delta \;=\;\frac{2\pi(2-2g)}{|V|}\;=\;\frac{\pi(1-g)}{2g}.
\]

Specifically, let $a=-\frac{\delta}{2}$
and define the base $8$ vertices
\begin{align*}
v_1 &= \left(\frac{1}{\tan{\left(\frac{\pi}{2g}\right)}},-1,1\right) \
v_2 = \left(\frac{1}{\tan{\left(\frac{\pi}{2g}\right)}},-1,-1\right) \
v_3 = \left(\frac{1}{\tan{\left(\frac{\pi}{2g}\right)}}-\tan{a},0,1\right) \\
v_4 &= \left(\frac{1}{\tan{\left(\frac{\pi}{2g}\right)}}-\tan{a},0,-1\right) \
v_5 = \left(\frac{1}{\tan{\left(\frac{\pi}{2g}\right)}}-\frac{1}{2}\sin{a}\tan{a},-1+\frac{1}{2}\sin{a},0\right) \\
v_6 &= \left(\frac{1}{\tan{\left(\frac{\pi}{2g}\right)}}-\frac{1}{2}\sin{a}\tan{a},1-\frac{1}{2}\sin{a},0\right) \
v_7 = \left(\frac{1}{\tan{\left(\frac{\pi}{2g}\right)}},1,1\right) \
v_8 = \left(\frac{1}{\tan{\left(\frac{\pi}{2g}\right)}},1,-1\right)
\end{align*}

For $k=0,\ldots,g-1$, all the vertices $v_{i,k}$ are obtained by rotating $v_i$ around the $z$–axis by $k\frac{2\pi}{g}$.

We define faces by
\begin{align*}
    &v_{1}v_{3}v_{7}v_{1,1}v_{3,1}v_{7,1}\cdots v_{1,g-1}v_{3,g-1}v_{7,g-1},\ v_{2}v_{4}v_{8}v_{2,1}v_{4,1}v_{8,1}\cdots v_{2,g-1}v_{4,g-1}v_{8,g-1},\\
    &v_{1,k}v_{3,k}v_{4,k}v_{5,k}v_{2,k},\ 
    v_{7,k}v_{3,k}v_{4,k}v_{8,k}v_{6,k},\ v_{1,k}v_{7,k}v_{6,k}v_{5,k},\ v_{2,k}v_{8,k}v_{6,k}v_{5,k},\\
    &v_{1,k}v_{2,k}v_{7,k}v_{8,k},\ v_{7,k}v_{8,k}v_{2,k+1}v_{1,k+1},\
    (k=0,\ldots,g-1)
\end{align*}
with the cyclic convention $v_{*,g}=v_{*,0}$.

\subsection{Orientable, without self-intersection, $6g$ vertices $(g\ge 4)$ }

Starting from the $V=8g$ arrangement, we compress the construction by switching to a regular $g$–gonal prism and opening a window on each of its lateral faces, so that each window still contributes one handle.
Each handle then contributes only $6$ vertices.
(Non–interference of adjacent windows requires $g\ge5$.)

\begin{figure}[htbp]
    \centering
    \includegraphics[width=0.2\linewidth]{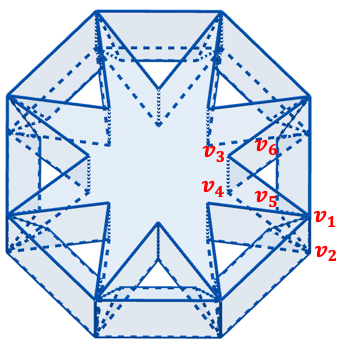}\hfill
    \includegraphics[width=0.2\linewidth]{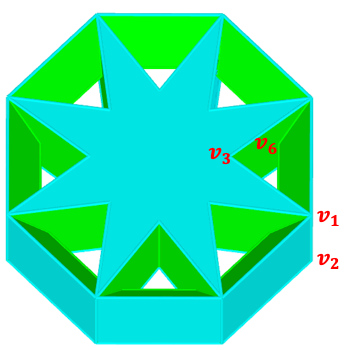}\hfill
    \includegraphics[width=0.2\linewidth]{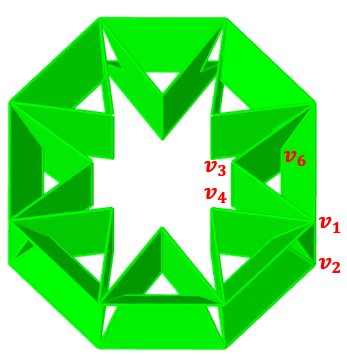}\hfill
    \includegraphics[width=0.8\linewidth]{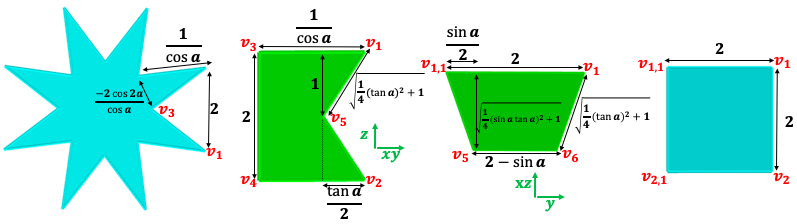}
    \caption{
        An embedded (without self-intersection)
     orientable CCP of genus $g=8$ with 48 vertices.
    }
    \label{fig:o_g8_v48}
\end{figure}

For every $g\ge5$ this yields an embedded orientable constant–defect polyhedral surface (CCP) with
\[
|V|=6g,\qquad |E|=\ 13g,\qquad |F|=\ 5g+2,\qquad 
\delta \;=\;\frac{2\pi(2-2g)}{|V|}\;=\;\frac{2\pi(1-g)}{3g}.
\]

Let $a=-\frac{\delta}{2}$ and
define the base $6$ vertices
\begin{align*}
v_1 &= \left(\frac{1}{\tan{\left(\frac{\pi}{g}\right)}},-1,1\right) &\quad
v_2 &= \left(\frac{1}{\tan{\left(\frac{\pi}{g}\right)}},-1,-1\right) \\
v_3 &= \left(\frac{1}{\tan{\left(\frac{\pi}{g}\right)}}-\tan{a},0,1\right) &
v_4 &= \left(\frac{1}{\tan{\left(\frac{\pi}{g}\right)}}-\tan{a},0,-1\right) \\
v_5 &= \left(\frac{1}{\tan{\left(\frac{\pi}{g}\right)}}-\frac{1}{2}\sin{a}\tan{a},-1+\frac{1}{2}\sin{a},0\right) &
v_6 &= \left(\frac{1}{\tan{\left(\frac{\pi}{g}\right)}}-\frac{1}{2}\sin{a}\tan{a},1-\frac{1}{2}\sin{a},0\right)
\end{align*}
For $k=0,\ldots,g-1$, all the vertices $v_{i,k}$ are obtained by rotating $v_i$ around the $z$–axis by $k\frac{2\pi}{g}$.

We define faces by
\begin{align*}
    &v_{1}v_{3}v_{1,1}v_{3,1}\cdots v_{1,g-1}v_{3,g-1},\ v_{2}v_{4}v_{2,1}v_{4,1}\cdots v_{2,g-1}v_{4,g-1},\\
    &v_{1,k}v_{3,k}v_{4,k}v_{2,k}v_{5,k},\ 
    v_{1,k}v_{6,k}v_{5,k}v_{1,k+1},\ 
    v_{2,k}v_{6,k}v_{5,k}v_{2,k+1},\ 
    v_{1,k}v_{2,k}v_{2,k+1}v_{1,k+1},\
    (k=0,\ldots,g-1)
\end{align*}
with the cyclic convention $v_{*,g}=v_{*,0}$.

\subsection{Orientable, without self-intersection, $7g-7$ vertices $(g=4,5,6)$}

Starting from the $V=6g$ scheme, we further reduce the vertex count by reusing inner window vertices to create an additional central cycle of windows, thereby decreasing the number of distinct vertices for small genera.

\begin{figure}[htbp]
    \centering
    \includegraphics[width=0.2\linewidth]{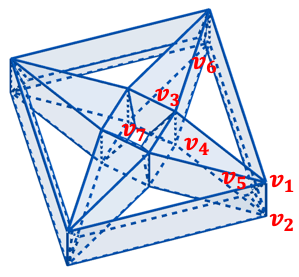}\hfill
    \includegraphics[width=0.2\linewidth]{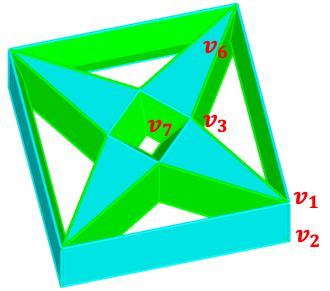}\hfill
    \includegraphics[width=0.2\linewidth]{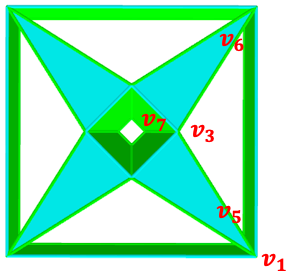}\hfill
    \includegraphics[width=0.2\linewidth]{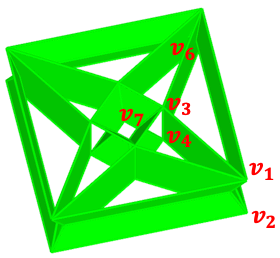}\hfill
    \includegraphics[width=0.8\linewidth]{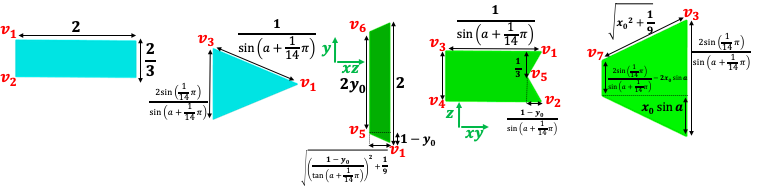}
    \caption{Top:
    An 
     orientable CCP without self-intersection of genus $g=5$ with 28 vertices,
     Bottom: Dimensions of the faces.
     A similar construction gives CCPs 
    for $g=4,5,6$.
    }
    \label{fig:o_g5_v28}
\end{figure}

For $g=4,5,6$ this yields embedded orientable CCPs with
\[
|V|=7g-7,\qquad |E|=18g-18,\qquad |F|=9g-9,\qquad 
\delta \;=\;\frac{2\pi(2-2g)}{|V|}\;=\;-\frac{4\pi}{7}.
\]

Choose $a=\dfrac{\pi}{g-1}$ so that $\delta=-2a$, and set
\[
x_0=\frac{\sqrt{2}\,\sin\!\bigl(\frac{\pi}{7}\bigr)}
           {3\,\sqrt{\cos\!\bigl(\frac{2\pi}{7}\bigr)-\cos(2a)}},
\qquad
y_0=1-\frac{\cos\!\bigl(\frac{3\pi}{14}\bigr)\,\sin\!\bigl(a+\frac{\pi}{14}\bigr)}
             {3\bigl(\sin\!\bigl(a+\frac{\pi}{14}\bigr)+\sin\!\bigl(\frac{3\pi}{14}\bigr)\bigr)}.
\]

Define the base $7$ vertices
\[
\begin{aligned}
v_1 &= \bigl(\frac1{\tan a},\,-1,\,\frac13 \bigr), 
&\quad 
v_2 &= \bigl(\frac1{\tan a},\,-1,\,-\frac13 \bigr),\\
v_3 &= \Bigl(\frac1{\tan a}-\frac1{\tan\bigl(a+\frac\pi{14}\bigr)},\,0,\,\frac13\Bigr),
&
v_4 &= \Bigl(\frac1{\tan a}-\frac1{\tan\bigl(a+\frac\pi{14}\bigr)},\,0,\,-\frac13\Bigr),\\
v_5 &= \Bigl(\frac1{\tan a}-\frac{1-y_0}{\tan\bigl(a+\frac\pi{14}\bigr)},\,-y_0,\,0\Bigr),
&
v_6 &= \Bigl(\frac1{\tan a}-\frac{1-y_0}{\tan\bigl(a+\frac\pi{14}\bigr)},\,y_0,\,0\Bigr),\\
v_7 &= \Bigl(\frac1{\tan a}-\frac1{\tan\bigl(a+\frac\pi{14}\bigr)}-x_0,\,0,\,0\Bigr),
\end{aligned}
\]
For $k=0,\ldots,g-2$, all the vertices $v_{i,k}$ are obtained by rotating $v_i$ around the $z$–axis by $k\frac{2\pi}{g}$.

The faces are
\begin{align*}
    &v_{1,k}v_{2,k}v_{1,k+1}v_{2,k+1},\ v_{1,k}v_{3,k}v_{3,k-1},\ v_{1,k}v_{5,k}v_{6,k}v_{1,k+1},\\ 
    &v_{1,k}v_{3,k}v_{4,k}v_{2,k}v_{5,k},\ v_{1,k+1}v_{3,k}v_{4,k}v_{2,k+1}v_{5,k+1},\ 
    v_{3,k}v_{7,k}v_{7,k+1}v_{3,k+1},\
    (k=0,\ldots,g-2)
\end{align*}
with the cyclic convention $v_{*,g-1}=v_{*,0}$.

\subsection{Non-orientable, odd genus}
\begin{figure}[htbp]
    \centering
    \includegraphics[width=0.3\linewidth]{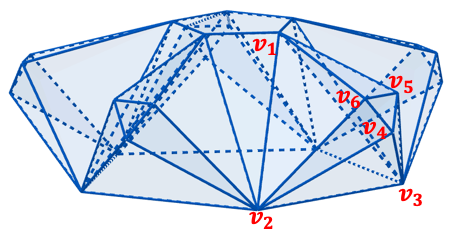}\hfill
    \includegraphics[width=0.3\linewidth]{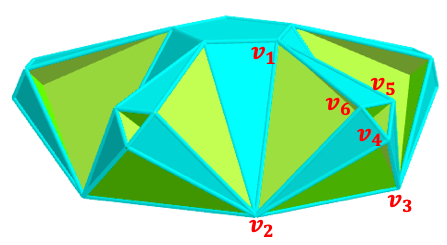}\hfill
    \includegraphics[width=0.2\linewidth]{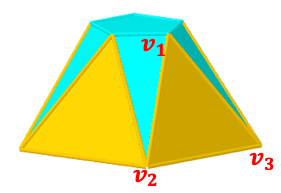}\hfill
    \includegraphics[width=0.5\linewidth]{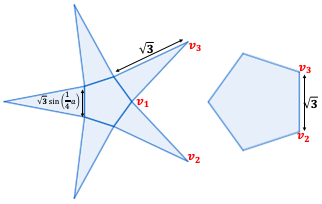}\hfill
    \includegraphics[width=0.3\linewidth]{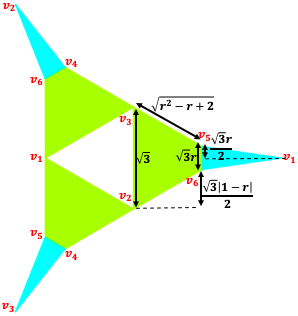}
    \caption{
    Top row left two: A non-orientable CCP of genus $g=5$ with 25 vertices,
    Top right: the base polyhedron,
    Bottom row left two: the net of of the base,
    Bottom-right: the net of the handle.
    The vertices $v_1v_2v_3$ are identified in the base and the handle.
    A similar construction gives CCPs 
    for odd $g\ge 3$.
    }
    \label{fig:n_g5_v25}
\end{figure}

For odd genera $3\le g\le11$, 
we construct a CCP (\Cref{fig:n_g5_v25} Top-left) with
\[
|V|=5g,\quad |E|=13g,\quad |F|=7g+2,\quad
\delta = \frac{2\pi(2-2g)}{|V|} = \frac{-2g+4}{5g}\,\pi.
\]
Taking a base polyhedron shown in the top-right of \Cref{fig:n_g5_v25}.
The base two vertices 
are given by
\[
\begin{aligned}
v_1 &= \Bigl(\frac{\sqrt3}{2}\tan\!\frac{5a}{4}
            -\sqrt{\frac94 - h_2^2}\,,\;0\,,\;h_2\Bigr),\\
v_2 &= \Bigl(\frac{\sqrt3}{2}\tan\!\frac{5a}{4}\,,\;
           -\frac{\sqrt3}{2}\,,\;0\Bigr),
\end{aligned}
\]
with
\[
h_2 = \frac{1}{2\cos\!\frac{5a}{4}}
      \sqrt{-3 + 6\cos a + 6\cos\frac{a}{2}
            -6\cos\frac{3a}{2} +6\cos\frac{5a}{2}},
            \quad
            a = \frac{2g-4}{5g}\,\pi.
\]
Note that the expression inside the radical becomes positive when $g\le 11$.
Other vertices are obtained by rotating the two about the $z$-axis by $k\frac{2\pi}{g}$, $k=1,\dots,g-1$. 

To the base, we attach $g$ handles $R(r,1)$ shown in the bottom-right of \Cref{fig:n_g5_v25}, where
\[
r \;=\;\frac12\Bigl(1-\sqrt{\frac{11+14\cos\!\frac{a}{2}}{1+2\cos\!\frac{a}{2}}}\;\tan\!\frac{a}{4}\Bigr).
\]

For genus $g\ge13$, we increase genus in steps of~2 by applying the drilling operation
to the $g=7$ construction, where 
we set $n=\frac{|V|}{\chi}=\frac{35}{5}=7$ by \eqref{eq:choice-of-n}.
% We verify that
% \[
% |V|=7g-14,\quad |E|=\frac{25g+7}{2},\quad |F|=\frac{9g+39}{2},\quad
% \delta=-\frac{2g-4}{7g-14}\,\pi.
% \]

\subsection{Non-orientable, even genus}

\begin{figure}[htbp]
    \centering
        \includegraphics[width=0.2\linewidth]{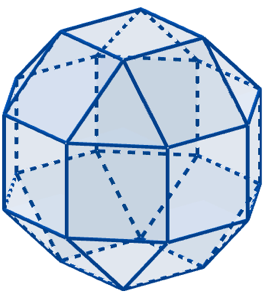}
    \includegraphics[width=0.2\linewidth]{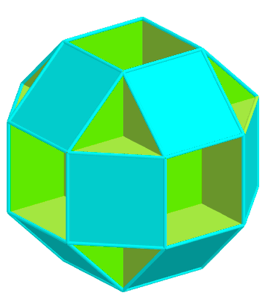}
    \caption{The rhombihexahedron, a non-orientable genus-$8$ CCP.
    CCPs with $4g - 8$ vertices for $g \geq 8$ 
can be obtained by drilling this with $\Pi(l,4)$.}
    \label{fig:n_g8_v24}
\end{figure}

\begin{figure}[htbp]
    \centering
    \includegraphics[width=0.2\linewidth]{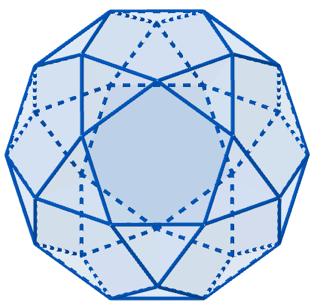}
    \includegraphics[width=0.2\linewidth]{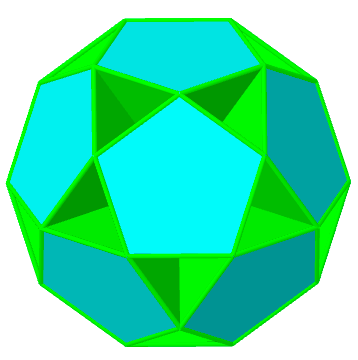}
    \caption{The small dodecahemidodecahedron, a non-orientable genus-$14$ CCP.}
    \label{fig:n_g14_v30}
\end{figure}

For genus $g = 2$, the CCP in \Cref{fig:n_g2_v9} has $9$ vertices.  
For genera $g = 4, 6$, the CCPs in \Cref{fig:n_g4_v12} have $6g - 12$ vertices.
For genera $g \geq 8$, CCPs with $4g - 8$ vertices
can be obtained by drilling with $\Pi(l,4)$ to \Cref{fig:n_g8_v24}.
However, as an exception, for genus $g = 14$, a CCP with fewer vertices, namely $30$, is known, as shown in \Cref{fig:n_g14_v30}.

%%%%%%%%%%%%%%%%%%%
\bibliographystyle{plain}
\bibliography{references}

\end{document}